\documentclass[11pt,notitlepage,twoside,a4paper]{amsart}
 \usepackage{amsfonts}
\usepackage[latin1]{inputenc}
\usepackage[cyr]{aeguill}
\usepackage{xspace}
\usepackage[polutonikogreek,french]{babel}

\usepackage{amsmath,amssymb,enumerate}

\usepackage{epsfig,fancyhdr,color}

\usepackage{amssymb}
\usepackage{amsmath,amsthm}
\usepackage{latexsym}
\usepackage{amscd}
\usepackage{psfrag}
\usepackage{graphicx}
\usepackage[latin1]{inputenc}
\usepackage[all]{xy}
\usepackage[mathcal]{eucal}

\definecolor{NoteColor}{rgb}{1,0,0}

\renewcommand{\textsc}{\textcolor{red}}
\newcommand*{\tg}[1]{\textgreek{#1}}

\newtheorem*{theorem 1}{\rm\bf Proposition 1}
\newtheorem*{theorem 2}{\rm\bf Proposition 2}

\theoremstyle{definition}

\theoremstyle{remark}

\def\interieur#1{\mathord{\mathop{\kern 0pt #1}\limits^\circ}}

\title[\`A partir du \emph{Timée} de Platon]{Mathématiques, musique et cosmologie : \`A partir du \emph{Timée} de Platon}
 
\author[Athanase Papadopoulos]{Athanase Papadopoulos}
\address{
Athanase Papadopoulos, Institut de Recherche Mathématique Avancée (CNRS et Université de Strasbourg), 7 rue René Descartes
67084 Strasbourg Cedex France
et GREAM (Groupe de Recherche Expérimental sur l'Acte Musical), Université de Strasbourg}
\email{papadop@math.unistra.fr}

\date{\today}

\begin{document}

   \maketitle

   \begin{abstract}
   \`A partir du Timée de Platon, nous présentons quelques réflexions sur la musique, la cosmologie et les mathématiques et leur influence mutuelle. 
   
   L'article est dédié au compositeur Walter Zimmermann.
   
   La version finale de l'article paraîtra dans le volume ``Les jeux subtils de la poétique, des nombres et de la philosophie. Autour de la musique de Walter Zimmermann", dir.P. Michel,  M. Andreatta et J.-L. Besada, Hermann, Paris, 2021.

   \bigskip
   
   \noindent {\sc Abstract.} 
   
   Based on Plato's Timaeus, we present some reflections on music, cosmology and mathematics and their mutual influence.

   The article is dedicated to the composer Walter Zimmermann.
   
   The final version of this article will appear in the volume ``Les jeux subtils de la poétique, des nombres et de la philosophie. Autour de la musique de Walter Zimmermann", ed. P. Michel, M. Andreatta and J.-L. Besada, Hermann, Paris, 2021.

   \bigskip
   
   AMS classification: 00A65, 00A30,  97M80, 01A20, 97A30

   \bigskip 
   
   Keywords: Mathematics and music, Plato, Timaeus
   
   \end{abstract}

      \section{Introduction}\label{s:Intro}
Je voudrais, dans les pages qui suivent, parler de musique chez Platon, en me plaçant du point de vue d'un mathématicien. Les idées que je présenterai me sont venues progressivement après que Pierre Michel m'ait demandé, en 2017, d'intervenir à un colloque consacré à l'\oe uvre du compositeur Walter Zimmermann. 
J'avais rencontré Zimmermann  trois ans auparavant et j'avais eu avec lui une conversation courte mais très dense. Ayant su que j'étais mathématicien, il m'avait parlé de l'une de ses compositions, me montrant une série de nombres sur laquelle il s'était appuyé en l'écrivant, et il avait mentionné Platon. Cette conversation m'avait particulièrement intéressé, car je voyais devant moi non seulement un musicien mais aussi un philosophe, et elle m'avait rappelé une phrase de Platon disant que le vrai musicien doit être philosophe. Aussi, quand Pierre Michel m'a demandé d'intervenir à ce colloque, j'ai pensé que ce serait  l'occasion pour moi de mettre en ordre quelques idées que je m'étais faites sur Platon à travers ses écrits, en particulier pendant les années où je m'y étais plongé de façon périodique pour préparer des cours, des conférences et des articles sur la relation entre les musique et les mathématiques et sur la philosophie des mathématiques.  Le colloque dédié à Zimmermann eut lieu le 27 mars 2018. Le présent article est fondé sur des notes que j'avais préparées pour la conférence que j'ai faite à ce colloque.

La série de nombres que Walter Zimmermann m'avait montrée est inspirée de celle du \emph{Timée} de Platon, c'est-à-dire celle que le Démiurge, dans ce dialogue, utilisa pour créer l'âme du monde, et c'est la même qui est employée  pour construire une gamme musicale, une longue gamme dont l'étendue est de quatre octaves et une sixte majeure. Le fait que la série de nombres qui a servi, dans le \emph{Timée}, à la création du monde soit la même que celle avec laquelle on obtient une gamme musicale est bien connu de tous ceux qui ont étudié ce dialogue. Cependant, il n'est pas possible de savoir si Platon a modelé la structure de l'âme du monde sur celle d'une gamme existante (une telle gamme n'apparaît dans aucun texte qui nous est parvenu avant le \emph{Timée}), et il n'est pas inconcevable de penser que c'est plutôt le développement ultérieur de la théorie de la musique, et en particulier celle de la construction des gammes, qui fut influencé par la construction par Platon (ou son Démiurge) de l'âme du monde.  

\begin{figure}[htbp]
\centering
\includegraphics[width=12cm]{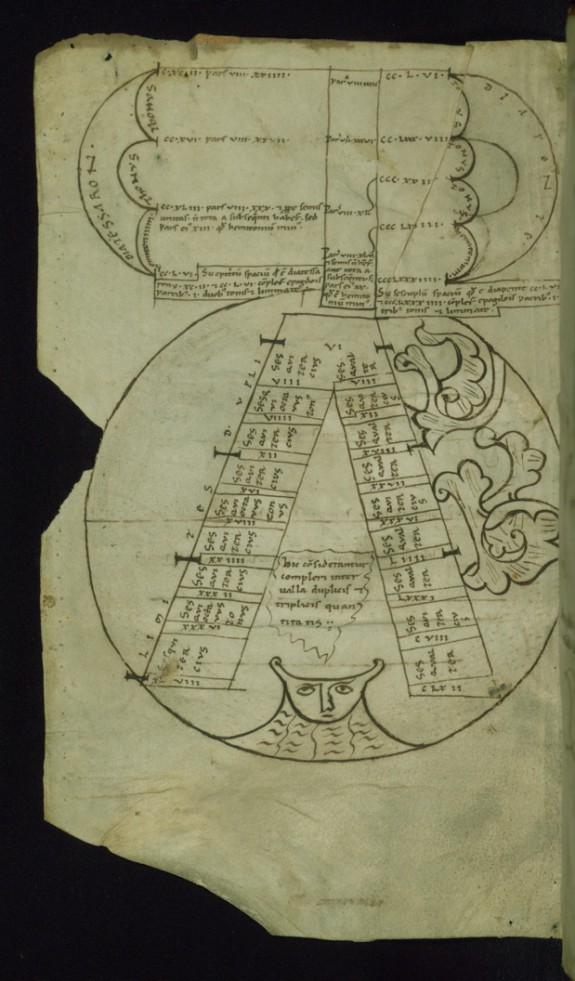}
\caption{\smaller Feuillet du \emph{Commentaire au songe de Scipion de Cicéron} de Macrobe, où l'on voit le fameux \emph{Lambda},  (Walters Art Museum, Baltimore).} \label{Macrobe}
\end{figure}

La série de Platon, ou les nombres qui l'engendrent, est souvent représentée sur les deux côtés d'une figure triangulaire ayant la forme de la lettre grecque $\Lambda$ (Lambda), une représentation qui date des  premiers commentateurs du \emph{Timée} et qui fut reprise par plusieurs autres  qui les suivirent. 
 Le feuillet reproduit en Figure \ref{Macrobe} ici contient une représentation de ce Lambda ; il est extrait d'un manuscrit du \emph{Commentaire au songe de Scipion de Cicéron} de Macrobe, philosophe et écrivain sné vers 370, l'un des derniers témoins de la fin de l'Antiquité romaine. Dans ce texte,  Scipion \'Emilien, jeune commandant de légion romaine, fait un rêve dans lequel il se voit expliquer la structure de l'univers par des a\"\i euls défunts lui parlant depuis des régions célestes où ils demeurent. Il y apprend en particulier que l'univers est composé de neuf sphères dont le mouvement produit une musique, appelée \emph{musique universelle} ou \emph{musique des sphères}, produite par des intervalles de notes dont les valeurs dépendent des distances entre les planètes et les autres astres connus. Le texte de Macrobe est fourni en références à Platon et à des philosophes platoniciens.

Dans les pages qui suivent, je voudrais présenter quelques idées de Platon sur la musique en mettant en valeur les relations qu'il a établies entre ce sujet, les mathématiques et la philosophie. J'inclurai la pensée de Platon en ce domaine, comme il est naturel de le faire, dans le cadre général de l'harmonie du monde, un thème qu'il développa dans le \emph{Timée} et qui fut poursuivi par d'autres après lui.

La suite de cet article est divisée en 5 parties, les sections \ref{s:math} à \ref{s:cercles}, dans lesquelles je développe ces thèmes, et elles sont suivies d'une conclusion. Le plan de ces sections est le suivant :

Dans la section \ref{s:math}, je mets en valeur deux textes pouvant nous éclairer sur a façon dont Platon concevait les mathématiques et la musique.

Dans la section \ref{s:Geometrie} j'évoque  le point de vue que Platon avait sur la géométrie. 

La section \ref{s:Musique} concerne la musique chez Platon. J'y parle de l'éducation musicale et du rôle de la musique comme moyen d'élévation pour l'âme. 

Dans la section \ref{s:Timee}, je me penche plus particulièrement sur les mathématiques et la musique dans le \emph{Timée}, l'écrit majeur de Platon dans lequel ce dernier traite de mathématiques et de musique de façon plus technique que dans la plupart de ses autres \oe uvres. C'est dans ce dialogue que l'on trouve une description détaillée des solides platoniciens et de leurs propriétés, et où appara\^\i t de façon explicite la gamme attribuée à Pythagore. Dans toute étude historique sur les polyèdres platoniciens ou sur les gammes grecques, une référence au \emph{Timée} s'impose. Le contexte de cette \oe uvre est cosmologique : Platon y décrit le Démiurge utilisant ces polyèdres et ces gammes dans la construction du corps et de l'âme de monde, ainsi que ceux de l'homme, mais il est surtout philosophique. Les tentatives d'interprétation  de cet écrit donnèrent lieu à des textes philosophiques fondamentaux répartis sur deux millénaires et demi ; j'en mentionne quelques-uns.

Dans la section \ref{s:cercles},  j'ai réuni quelques passages et quelques images pour illustrer la manière dont le cercle, comme figure géométrique, est utilisé, depuis plus de deux millénaires, par des mathématiciens-musiciens pour représenter des concepts géométriques, musicaux et cosmologiques.

\section{Platon, les mathématiques et la musique} \label{s:math}

Platon était en premier lieu mathématicien, et le regard qu'il portait sur le monde était celui d'un mathématicien. Je voudrais commencer par citer Aristoxène de Tarente, le grand théoricien de la musique du IV$^{{\mathrm{e}}}$ siècle avant J.-C. 
 qui a commenté cet aspect de Platon. Aristoxène connaissait les mathématiques platoniciennes et ils savaient comment fonctionne l'esprit d'un mathématicien. J'ai utilisé l'expression \og mathématiques platoniciennes\fg, en considérant qu'elle a un sens, mais j'aurais pu dire mathématiques tout court car l'expression désigne une tradition qui se perpétue encore aujourd'hui. Plus loin dans l'article, je ferai des références plus directes à des textes de Platon lui-même, pour décrire ce genre de mathématiques.

Dans un passage dans lequel il établit une comparaison entre le maître de l'Académie et celui du Lycée,  Aristoxène écrit dans son \emph{Traité d'harmonie} : \og Platon all\'echait son
auditoire en lui promettant une conf\'erence sur le Bien, pour finalement ne
l'entretenir que d'Arithm\'etique, de G\'eom\'etrie, d'Astronomie, et conclure
enfin que ``le Bien est Un"\fg \cite[p. 98]{A-Belis}. Ce fait, qui pourrait paraître à première vue anecdotique, nous éclaire sur la manière dont Platon expliquait sa vision du monde,  une vision selon laquelle seul le discours mathématique peut exprimer des idées intelligibles, c'est-à-dire que l'on peut atteindre par l'intelligence et non par les sens. Pour lui, des notions comme celles du Bien et du Beau ne peuvent être liées à des objets physiques parce que, disait-il, tout ce qui est physique est périssable et incompatible avec l'idée du Bien et du Beau.

 Les écrits de Platon sont parsemés de références aux mathématiques et à la musique et aux relations entre ces deux domaines. 
       Dans le \emph{Philèbe} (17d-e), Platon fait dire à Socrate que c'est par la connaissance des intervalles de musique, par celle des systèmes qui en résultent, par l'étude de la classification que les Anciens en ont fait ainsi que par la connaissance du rythme et de la mesure que l'on acquiert l'intelligence et que l'on devient savant.
\`A ce propos, je voudrais citer Théon de Smyrne qui consacre une partie importante de son \emph{Exposition des connaissances math\'ematiques utiles pour la lecture de Platon}    \cite{Theon}  à la musique. Dans l'introduction, il rappelle que Platon décrit le musicien comme un penseur qui est en même temps philosophe, insistant sur les qualités intellectuelles de tempérance et générosité que doit avoir celui qui dans la cité se destine à la musique. Dans le même passage, les qualités intellectuelles et les vertus morales de l'homme en général sont décrites en termes musicaux : harmonie, accord, cadence dans le rythme, etc.   Théon écrit: \og [Platon] montre que le seul philosophe est réellement musicien, tandis que celui qui est vicieux et méchant est étranger aux Muses. Car, dit-il, la vraie et sincère probité des m\oe urs, cette vertu qui consiste dans le bon et honnête règlement de notre vie, suit la droite raison, c'est-à-dire l'usage conforme à la raison. Il ajoute que les compagnons de la droite raison sont la décence, la cadence et l'accord, la décence dans le chant, l'accord dans l'harmonie, la cadence dans le rythme. Par contre, l'improbité ou la corruption des m\oe urs est essentiellement liée à la perversion de la raison, c'est-à-dire à l'usage corrompu de la raison, et ses compagnons sont l'indécence, la confusion et le désaccord dans tout ce qu'on fait, de soi-même ou par imitation, de sorte que celui-là seul est musicien qui a de bonnes m\oe urs et, comme on le voit par ce qui précède, il est aussi le vrai philosophe, si toutefois, dès les premières années de son adolescence, quand on lui eut appris la musique, il prit des habitudes de décence et d'ordre, car la musique joint un plaisir innocent à l'utilité. Il est impossible, dit Platon, que celui-là devienne musicien parfait, qui n'a pas en tout des habitudes de bonne éducation, qui n'a pas les idées de décence, de noblesse d'âme et de tempérance. Il doit reconnaître que ces idées se retrouvent partout et ne les mépriser ni dans les petites choses ni dans les grandes. Car c'est au philosophe qu'il appartient de connaître les idées, et personne ne connaîtra la modestie, la tempérance et la décence, s'il est lui-même immodeste et intempérant. Mais les choses qui font l'ornement de la vie humaine, le beau, l'harmonieux, l'honnête, tout cela est l'image de cette beauté, de cet accord, de ce bel ordre éternel et qui a une existence véritable, c'est-à-dire que ces choses sensibles sont les caractères et l'expression des choses intelligibles ou des idées. \fg

   \section{Géométrie} \label{s:Geometrie}

    Platon n'aimait pas le sens étymologique du mot géométrie, comme \og mesure de la terre\fg. Il écrit dans  l'\emph{\'Epinomis} 990d \cite[p. 1160]{Epinomis} que ce nom est \og parfaitement ridicule\fg.
 Pour lui, une telle interprétation de cette science nous plonge dans des réalités pratiques, alors que la géométrie devrait être une voie vers la connaissance du monde des idées, nous donnant une sensation des réalités immuables. La géométrie nous fait pénétrer le monde des formes, comme opposé à celui de la matière.  Cela nous renvoie aussi au passage d'Aristoxène que nous avons mentionné : ayant annoncé des discours sur le Bien et sur le Beau, Platon offrait finalement à ses auditeurs des exposés de mathématiques qui constituent eux-mêmes une méditation sur le Bien et sur le Beau. Dans la \emph{République} 526e \cite[p. 118]{Republique}, parlant des études qui devraient être à la base d'une bonne éducation, Platon fait dire à l'un de ses personnages que \og cette étude [de la géométrie] vise en quelque chose à ce but élevé, qui est de faire que nous voyions plus aisément la nature du Bien ; but auquel, disons-nous, vise toute étude qui force l'âme à se détourner dans la direction de cette région sublime où réside ce qui, dans le réel, possède la plus haute béatitude et dont il faut, à tout prix, que l'âme ait eu la vision. \fg  Dans le même passage, il déclare : \og Que la géométrie force à contempler la réalité, elle nous intéresse ; qu'elle force à contempler la génération, elle ne nous intéresse pas\fg \cite[p. 1119]{Republique}. Pour Platon, la réalité, c'est le monde des idées, et la génération, c'est le monde sensible, celui des choses qui naissent et qui ensuite disparaissent. Toujours d'après Platon, la géométrie, comme connaissance, ne s'occupe pas du monde physique. Même s'il est vrai que ce monde est une imitation du monde immuable des idées, il ne peut être l'objet d'une connaissance car il est en perpétuel changement. Dans la suite du passage que nous venons de citer, Platon écrit que \og la nature de cette science s'oppose totalement au langage de certains qui la manient\fg , un langage, dit-il \og tout à fait risible et qui sent la servilité\fg, celui \og de gens qui pratiquent une action et dont la pratique est le but\fg, ceux qui parlent de travaux pratiques, comme le fait de \og carrer\fg, de \og tendre le long de\fg, de \og poser en plus de \fg, alors que l'objet de cette science est la connaissance, une connaissance \og de ce qui toujours existe, et non pas de ce qui, à un moment donné, commence ou finit d'être quelque chose\fg.

Une question de géométrie qui  vient naturellement à l'esprit dans ce contexte platonicien est celle de la classification des polyèdres convexes réguliers, que l'on appelle aussi polyèdres platoniciens (les deux terminologies sont courantes). Je voudrais d'abord rappeler en quelques mots de quoi il s'agit. 

Un  polyèdre platonicien est une surface dans l'espace à trois dimensions obtenue en recollant le long de leurs bords des figures planes bordées par des polygones réguliers (des triangles équilatéraux, des carrés, des pentagones réguliers, des hexagones réguliers, etc.) identiques entre eux. Ces figures deviennent les faces du polyèdre. De plus, on demande que ces faces soient recollées entre elles de telle sorte qu'il y ait toujours le même nombre de faces qui touchent chaque sommet du polyèdre. Cette dernière propriété fait partie de la \emph{régularité} du polyèdre : non seulement les faces sont toutes identiques entre elles, mais aussi, les sommets sont identiques entre eux. On demande de cette façon une symétrie maximale. Enfin, un polyèdre platonicien est convexe, c'est-à-dire qu'il est totalement situé d'un seul côté par rapport à n'importe quel plan contenant une de ses faces. On montre que sous ces conditions, il n'y a que cinq polyèdres. Ainsi, très peu de pentagones réguliers peuvent être utilisés pour fabriquer des polyèdres convexes réguliers : il n'y a que le triangle, le carré et le pentagone.

Classer les polyèdres platoniciens, savoir qu'il y en a exactement cinq, est un théorème non trivial, que Platon connaissait et que les Pythagoriciens avaient découvert avant lui\footnote{Heath, dans son édition des \emph{\'Eléments}, au début du livre XIII, a écrit des notes historiques intéressantes sur la construction des polyèdres platoniciens. Un historien des sciences a publié un article de 78 pages sur les mathématiques dans le \emph{Timée} \cite{Vitrac}. Il écrit (p. 22) : \og Bien entendu, il n'est pas question de suggérer que des méthodes mobilisant dans leur mise en \oe uvre effective des outils mathématiques sophistiqués comme la théorie des équations différentielles, des méthodes de résolution par approximation et des moyens de calcul très puissants, aient existé
à l'époque de Platon\fg. L'expression \og outils mathématiques sophistiqués\fg renvoie ici au contexte de quelqu'un qui est peu au courant de ce que sont réellement les mathématiques et qui n'en connaît que ce que les auteurs en sciences sociales appellent \og outils mathématiques\fg ; en réalité, ces outils que sont \og équations différentielles\fg, les \og méthodes de résolution par approximation\fg et les\og moyens très puissants\fg auxquels il fait référence sont, du point de vue conceptuel, triviaux comparés aux idées géométriques profondes sur la classification et la construction des polyèdres réguliers que l'on trouve dans le \emph{Timée}. Je ne crois pas que le but de Platon était de faire de la \og modélisation mathématique\fg comme le pense cet auteur ; il faisait tout simplement un exposé de mathématiques, comme il avait l'habitude de le faire quand il parlait de philosophie.}. Les preuves modernes de ce théorème font appel à la formule de Descartes-Euler\footnote{C'est la formule qui dit que pour n'importe quelle décomposition polygonale de la sphère, le nombre de faces moins le nombre d'arêtes plus le nombre de sommets est une constante, ne dépendant pas de la décomposition. La formule fut longtemps attribuée à Euler, mais on sait maintenant que Descartes la connaissait avant lui.}, c'est-à-dire à des mathématiques du dix-septième ou du dix-huitième siècles. Si l'on ne connaît pas cette formule, il faut faire, pour conduire la classification, des considérations combinatoires qui ne sont pas faciles (qu'Euclide connaissait). Il faut aussi utiliser les symétries de la sphère, un fait sur lequel Platon insistait. Cette classification est contenue dans le dernier livre des \emph{\'Eléments}, et en quelque sorte elle en est le couronnement. 

 On ne connaît pas le nom de celui qui fut le premier à avoir dégagé la notion de polyèdre platonicien, encore moins celui qui en découvrit les premières propriétés importantes. Une des raisons de ce manque d'information est que les Pythagoriciens constituaient une confrérie où l'on mettait tout en commun :  pas seulement les biens matériels, mais aussi les découvertes mathématiques. Les théorèmes obtenus étaient attribués à la communauté dans son ensemble. Ces découvertes revêtaient une importance particulière pour eux, à tel point qu'il fallait les garder secrètes. Jamblique, l'écrivain néo-pythagoricien syrien qui vécut aux II$^{{\mathrm{e}}}$-III$^{{\mathrm{e}}}$ siècles après J.-C., raconte dans sa \emph{Vie de Pythagore} qu'Hippase le Pythagoricien, \og parce qu'il avait été le premier à divulguer par écrit comment on pouvait construire une sphère à partir de douze pentagones, périt en mer pour avoir commis un acte d'impiété \fg \cite[p. 51]{Jamblique}. Le même auteur rapporte que \og certains ont dit que celui aussi qui avait révélé ce qui concerne l'irrationalité et l'incommensurabilité a subi le même sort \fg \cite[p. 132]{Jamblique}. Que ces histoires aient ou non une base historique quelconque n'est pas essentiel pour nous, et d'ailleurs on ne pourra jamais les confirmer ni les infirmer ; mais le fait qu'elles soient passées dans la littérature témoigne de l'importance que revêtaient les polyèdres réguliers pour les Pythagoriciens et pour Platon qui en était l'héritier. Ce dernier, dans le \emph{Timée}, nous offre un long passage dans lequel il décrit la construction de ces polyèdres ; nous en reparlerons dans la section \ref{s:Timee}.

\begin{figure}[htbp]
\centering
\includegraphics[width=12.5cm]{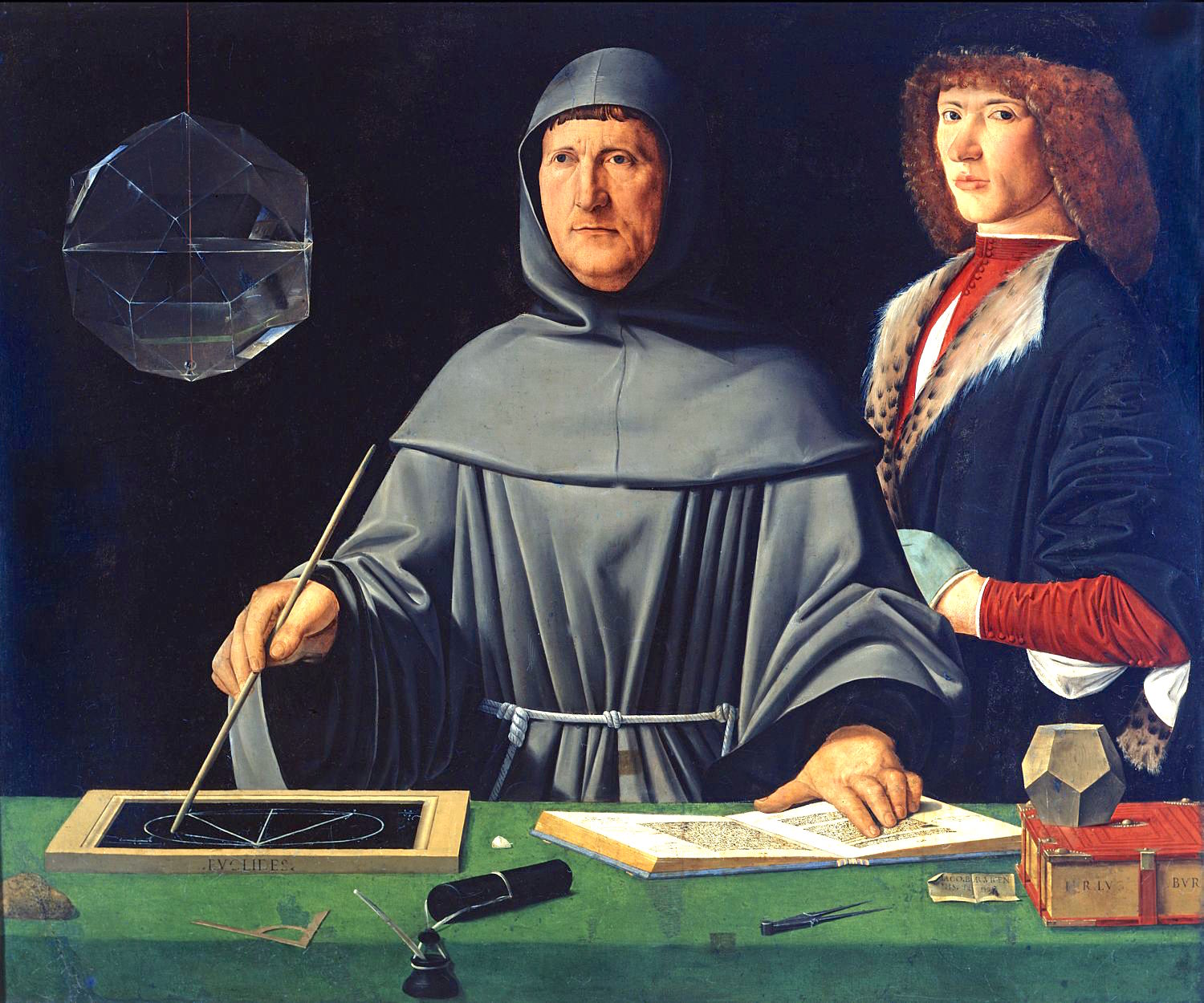}
\caption{\smaller Portrait de Luca Pacioli avec un étudiant (peut-être Albrecht Dürer) avec, sur la table, un dodécaèdre régulier, et accroché en l'air, un polyèdre semi-régulier formé de triangles et de carrés. Peinture attribuée à Jacopo de' Barbari  (Musée de Capodimonte, Naples).} \label{Barbari}
\end{figure}

\begin{figure}[htbp]
\centering
\includegraphics[width=12.4cm]{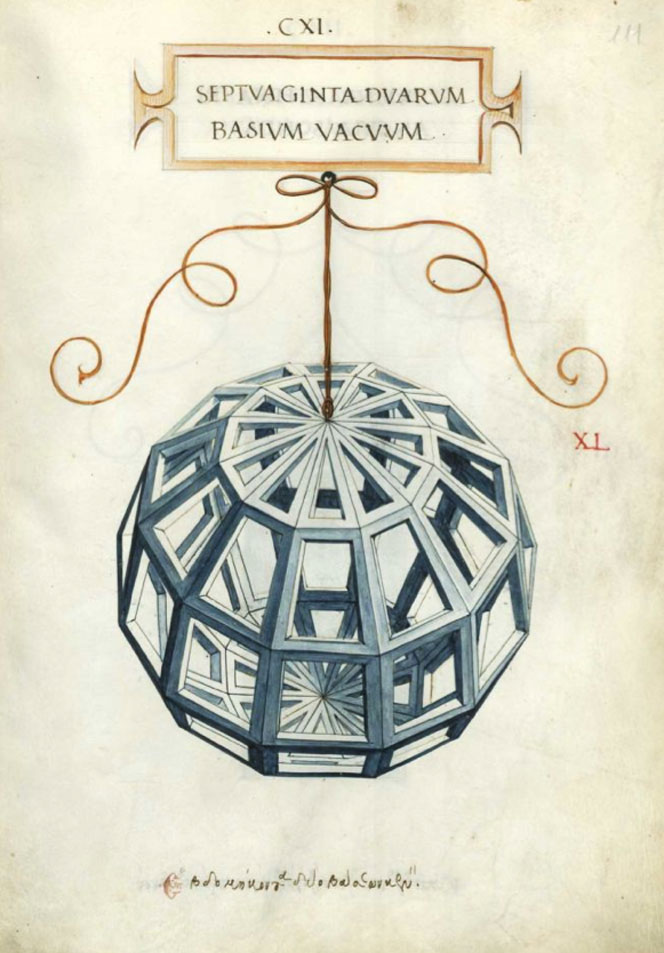}
\caption{\smaller Dessin tiré du \emph{De divina proportione} de Luca Pacioli (Venise, 1469)} \label{Pacioli}
\end{figure}

\begin{figure}[htbp]
\centering
\includegraphics[width=12.5cm]{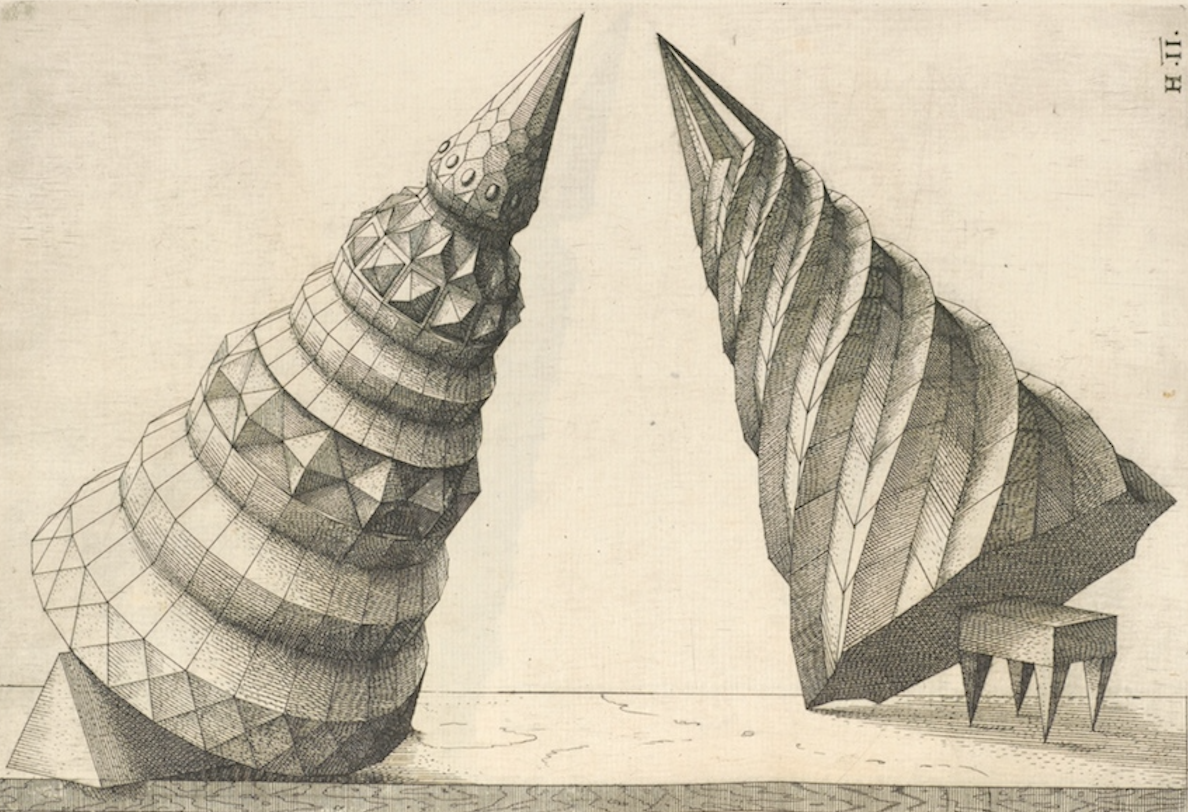}
\caption{\smaller Polyèdres. Dessin tiré du livre \emph{Perspectiva corporum regularium} 
	de  Wenzel Jamnitzer  (Nürnberg, 1568). Bibliothèque du Getty Research Center, Los Angeles.} \label{Jamnitzer}
\end{figure}

 Le mathématicien, quand il a  résolu une question, en invente une nouvelle, parfois plus générale ou plus compliquée.
Ainsi, des polyèdres convexes réguliers, on passe naturellement à ceux qui sont semi-réguliers, c'est-à-dire ceux dont les faces sont faites de polygones réguliers mais pas forcément identiques les uns aux autres. On demande, de plus, comme pour les polyèdres réguliers, une symétrie par rapport des sommets. Les faces, par exemple, pourraient être un mélange de triangles équilatéraux et de carrés, ou de carrés et de pentagones réguliers, etc. C'est à Archimède, qui vécut environ un siècle et demi après Platon, que l'on doit la classification de ces polyèdres. Si l'on exclut les polyèdres convexes réguliers,  il y en a essentiellement treize\footnote{L'ouvrage d'Archimède qui contient ce résultat est perdu, mais Pappus en fait rapport dans sa \emph{Collection}, un ouvrage extrêmement riche en résultats mathématiques et qui est en même temps une source fondamentale pour la connaissance de l'histoire des mathématiques grecques. La classification de ces polyèdres, que l'on appelle archimédiens, est contenue dans le livre V de la \emph{Collection}, voir \cite{Pappus}, t. 1,  p. 273 ss.}. Les peintres de la Renaissance comme D\"urer, Leonard de Vinci, Piero della Francesca et Luca Pacioli, ont laissé de très beaux dessins et tableaux de polyèdres, réguliers, semi-réguliers et d'autres encore plus généraux, et ils en ont utilisé certains comme modèles pour la construction de dômes et de coupoles. La figure \ref{Barbari} reproduit un tableau de Jacopo de' Barbari, représentant, au centre, Luca Pacioli et où l'on voit, posé sur la table, un dodécaèdre régulier, et accroché en l'air, un polyèdre archimédien dont les faces sont constituées par des triangles équilatéraux et des  carrés. La figure \ref{Pacioli} représente un polyèdre  tiré du livre \emph{De divina proportione} de Luca Pacioli. C'est un polyèdre non régulier, dessiné pour son esthétique. Enfin, la figure \ref{Jamnitzer} provient de l'\oe uvre du dessinateur, graveur et orfèvre du XVI$^{\mathrm{e}}$ siècle  Wenzel Jamnitzer.
 
 Aujourd'hui les polyèdres jouent un rôle important dans la conception des salles de concert, mais il ne semble pas que le sujet ait été très exploité du point de vue de la composition musicale, même si beaucoup d'idées sur ce thème peuvent venir naturellement à l'esprit (on pense quand même aux \emph{Polytopes} de Xenakis).
  
Signalons pour terminer cette section que les mathématiques platoniciennes, ce n'est pas seulement la géométrie, mais c'est aussi la théorie des nombres. Dans l'\emph{\'Epinomis}  976c-e, il écrit, en parlant de ce domaine :  \og C'est l'unique connaissance qui, si elle désertait notre humaine nature ou si elle ne lui était pas présente, ferait de la race des hommes la plus déraisonnable et la plus insensée des espèces animales actuellement existantes\fg \cite[p. 1139]{Epinomis}. La musique, dans l'Antiquité grecque, faisait partie de la théorie des nombres. Dans la prochaine section, on parlera donc de musique.
  
\section{Musique}\label{s:Musique}

De même qu'il l'était pour les mathématiques, Platon, concernant la musique, était l'héritier des Pythagoriciens.
   La musique joue un rôle prépondérant dans ses écrits, et pour lui, le sujet avait, comme pour les Pythagoriciens, une dimension éthique, esthétique, physique et cosmique. En effet, la musique est l'une des manifestations de  l'harmonie dans l'univers, le mot harmonie évoquant ici une notion d'ordre, désignant une combinaison réussie, ou une construction coordonnée\footnote{La notion de l'harmonie dans l'univers, qui est un thème fondamental dans les écrits de Platon, fut largement explorée dans la  littérature qui a suivi. On pense à Cicéron :  \og Mérite-t-il vraiment le nom d'homme, celui qui, en présence de tant de mouvements bien réglés, d'un ordre si parfait régnant au ciel, des liens unissant de façon si harmonieuse toutes les parties du monde les unes aux autres, se refuse à croire à une raison ordonnatrice, prétend mettre au compte du hasard un arrangement calculé de façon si savante que notre science en est déconcertée? \fg (\emph{De la nature des dieux}, Livre II). On pense aussi à Dante : \og Quand la roue [le firmament] qu'anime son désir éternel pour toi [pour Dieu] attira mon attention par son harmonie que tu tempères et répartis [\ldots]\fg (\emph{La Divine Comédie}, Le Paradis, Chant I). On pense aussi à Shakespeare qui, dans l'acte V, scène 1 du \emph{Marchand de Venise}, fait résonner une musique des sphères dans le but de faire sentir l'harmonie entre le macrocosme et le microcosme, un thème cher à Platon et aux Pythagoriciens. On pourrait multiplier les exemples.}. C'est pour cette raison que le thème de la musique est sous-jacent de façon fondamentale au \emph{Timée}, le dialogue dans lequel Platon décrit la construction du monde, et c'est pour cela que le mot harmonie y est souvent cité. Théon écrit dans son \emph{Exposition} : \og  Les Pythagoriciens dont Platon adopte souvent les sentiments, définissent aussi la musique comme une union parfaite totale de choses contraires, l'unité dans la multiplicité, enfin l'accord dans la discordance. Car la musique ne coordonne pas seulement le rythme et la modulation, elle met l'ordre dans tout le système ; sa fin est d'unir et de coordonner, et Dieu aussi est l'ordonnateur des choses discordantes, et sa plus grande \oe uvre est de concilier entre elles, par les lois de la musique et de la médecine, les choses qui sont ennemies les unes des autres. C'est aussi par la musique que l'harmonie des choses et le gouvernement de l'univers se maintiennent ; car ce que l'harmonie est dans le monde, la bonne législation l'est dans l'\'Etat, et la tempérance l'est dans la famille. Elle a, en effet, la puissance de mettre l'ordre et l'union dans la multitude\fg \cite[p. 19]{Theon}.

  On dit que Pythagore était capable écouter l'harmonie universelle, une musique produite par le mouvements des astres, et ses biographes ont ajouté que si cette activité nous est étrangère, à nous humains ordinaires,  c'est à cause de l'étroitesse de nos facultés\footnote{Cf. par exemple ce texte du philosophe néo-platonicien Prophyre de Tyr (III$^{{\mathrm{e}}}$ siècle après J.-C.), parlant de Pythagore : \og Il écoutait l'harmonie de l'univers, sensible qu'il était à l'harmonie universelle des sphères et des astres qui s'y meuvent, que la médiocrité de notre nature nous empêche d'entendre\fg \cite[p. 30]{Porphyre}.}. Il est probable que cela soit une légende, mais il est possible aussi que cette affirmation soit une façon déguisée --- on rappelle que l'enseignement des Pythagoriciens était ésotérique --- de dire que le mouvement des astres est périodique, et que comme tout mouvement périodique, il a une fréquence, que toute fréquence correspond à une certaine note musicale, même si cette note ne peut pas être entendue, et que Pythagore comprenait cela. On sait par ailleurs que les analogies entre la musique et l'astronomie faisaient partie des sujets favoris des Pythagoriciens, et que les deux domaines se sont mutuellement influencés pendant plus de deux mille ans. Par exemple, ayant constaté que le nombre des grands luminaires célestes (la Lune, le Soleil et les cinq planètes connues) était de sept, les Pythagoriciens conçurent une lyre ayant exactement sept cordes\footnote{Cet instrument, connu du temps de Philolaos,  était accordé suivant une gamme formée de deux tétracordes successifs ayant une note commune.   Nicomaque,  Pline, Ptolémée, Boèce et d'autres ont décrit les gammes construites à partir de tels tétracordes ; voir à ce sujet l'article intéressant de Reinach \cite{Reinach}.}.
  
Inversement, les théories sur l'harmonie du monde jouèrent un rôle certain dans la composition musicale. Un grand nombre de pièces inspirées par l'harmonie céleste, et d'autres ayant pour but de faire sentir que l'univers dont nous faisons partie est harmonieux virent le jour pendant la Renaissance. L'expression \emph{harmonia mundi} date de cette époque.  Il faut ajouter que le thème de l'harmonie céleste fut christianisé, et il est difficile parfois de  séparer les musiques qui en sont inspirées de celles que nous appelons \og sacrées\fg, avec une connotation chrétienne. 
   Parmi les compositeurs qui ont exploité le thème de l'harmonie de l'univers, on peut citer Thomas Tallis, William Byrd, Thomas Morley, Michel Maier  et Carlo Gesualdo.
    Nous avons déjà mentionné l'harmonie des sphères dans l'acte V du \emph{Marchand de Venise} de Shakespeare. On pense aussi à l'opéra \emph{Les Planètes} de Lully,  ainsi qu'au \emph{Mondo della luna} et à la \emph{Création} de Haydn. Le thème de la musique céleste fut ravivé par les compositeurs au vingtième siècle. On pense à certaines compositions de Schoenberg, Gustav Holst, Philip Glass, Stockhausen, et, bien sûr, Walter Zimmermann. 
Dans un texte qui accompagne la pièce pour piano de ce dernier,  \emph{Wüstenwanderung} (Voyage dans le désert), composée en 1986, on lit que cette \oe uvre \og décrit la création de l'âme du monde d'après le \emph{Timée} de Platon, devenant de plus en plus compliquée, pour s'effondrer sous sa propre complexité, qui est devenue celle d'une machine.\fg

  Platon  développa le thème de la relation entre l'harmonie musicale et l'harmonie céleste dans plusieurs écrits et il y intégra l'harmonie de l'âme, une notion liée, selon lui, aux qualités de raison et de sagesse.  Dans les  \emph{Lois} 689d, on lit  : \og Parmi les harmonies concertantes, la plus belle et la plus haute est celle qui consiste en la plus haute sagesse, celle à laquelle participe l'homme qui vit d'une manière raisonnable \fg  \cite[p. 725]{Lois}. 
 Dans la 
  \emph{République} 443d-e, il parle de réaliser entre les trois parties de l'âme \og un accord, ni plus ni moins que s'il s'agissait des termes d'un accord musical, entre haute, basse et moyenne, sans compter tels autres termes qu'on peut introduire entre ceux-là ; opérant la liaison de tout cela et, avec une multiplicité, nous faisant unité, tempérament, harmonisé\fg. Il y aurait beaucoup d'autres exemples à citer.

  Il faut rappeler ici que pour Platon, le mot \og musique\fg désignait essentiellement la musique théorique. 
De même qu'il critiquait les géomètres qui avaient l'esprit dirigé vers les choses pratiques, il critiquait ceux parmi les musiciens qui, disait-il, passaient leur temps, l'oreille tendue, à écouter la justesse de tel ou tel intervalle. Il opposait le musicien-philosophe, qu'il considérait comme le vrai musicien,  à l'artiste-sophiste\footnote{Dans la bouche de Platon, un sophiste est un imposteur.}, un virtuose de l'illusion et des faux-semblants, spécialiste des apparences. Alors que le premier attire l'âme vers le Bien, le second, pris par des tâches musicales empiriques et non raisonnées, la détourne vers ce qui est une chimère. Il est \og en possession d'un certain art de l'apparence illusoire\fg\footnote{\emph{Sophiste} 239d \cite[p. 292]{Sophiste}.} et de celui de \og fabriquer des simulacres\fg\footnote{\emph{Sophiste} 260d \cite[p. 325]{Sophiste}.}. Pour Platon, la musique, et plus généralement l'art, n'est pas une technologie, mais comme les mathématiques, c'est une méditation ; elle attire l'âme vers ce qui est beau et, par conséquent, vers ce qui est immuable et sans changement, car ce qui est beau ne peut être qu'éternel\footnote{\emph{Timée} 28b \cite[p. 444]{Timee}.}. 
On peut citer ici un passage de la \emph{République} 531a, où l'on trouve une critique des gens qui s'occupent de mesurer des consonances, savoir si elles sont sensibles ou pas à l'ouïe ; ces gens ressemblent, nous dit Platon, aux astronomes qui passent leur temps à faire des mesures mais dont le travail n'aboutit pas  : \og  --- Par tous les Dieux! fit-il, risiblement : donnant des noms distincts à certaines concrétions des sons, tendant l'oreille comme s'ils étaient à l'affut de ce qui se dit chez les voisins, les uns prétendent discerner encore par l'ouïe une différence, moyenne entre deux tons, et qui est le plus petit intervalle dont il faille se servir pour mesurer ; les autres le contestent et prétendent qu'il est pareil aux sons antérieurement émis  : les uns comme les autres donnant aux oreilles la prééminence sur l'intellect. --- Tu penses, toi, repris-je, à ces excellentes gens qui donnent aux cordes mille tracas, qui les mettent à la torture en les tordant au moyen des clefs! \fg  \cite[p. 1124]{Republique}.

La musique, pour Platon, sert à élever l'âme et non pas à  la rendre inquiète. Comme les autres activités artistiques, elle ne doit pas être une fin en soi ; elle est pratiquée pour les vertus qu'elle suscite dans l'âme. Il insiste naturellement sur la musique comme élément essentiel de l'éducation des jeunes.  Son livre des \emph{Lois} (810a-813) régit strictement l'apprentissage de la musique. On y  lit par exemple : \og Quand [l'enfant] aura treize ans, ce sera le bon moment pour commencer d'aborder l'étude de la lyre, étude à laquelle il passera trois autres années. Défense au père, aussi bien qu'à l'enfant lui-même et que ces études lui plaisent ou qu'il les déteste, d'augmenter ou d'abréger cette durée, de consacrer à cela un temps illégal, en plus ou en moins\fg \cite[p. 895ss]{Lois}. Il y aurait plusieurs autres textes à citer sur ce sujet. Dans la \emph{République} 401d, il écrit que
\og le motif pour lequel la culture musicale est d'une excellence souveraine est que rien ne plonge plus profondément au c\oe ur de l'\^ame que le rythme et l'harmonie ; que rien ne la touche avec plus de force en y portant l'harmonieuse élégance qui en fait la noblesse, dans le cas où cette culture a été correctement conduite
\fg \cite[p. 957]{Republique}. Dans l'\emph{\'Epinomis}  989b, il écrit :  \og [\ldots]  C'est qu'une [âme] qui reçoit, et avec gentillesse, une juste proportion de lenteur et de la qualité opposée, sera sans doute d'humeur facile, aimera la vaillance autant qu'elle sera docile aux conseils de la modération et, ce qu'il y a de capital, quand à de semblables naturels appartient une âme capable d'apprendre et de se souvenir, cette âme sera capable de trouver un vif plaisir aux questions dont il s'agit, au point de devenir fervente amie de l'étude \fg  \cite[p. 1158]{Epinomis}.

  La musique, pour Platon, est un don des dieux. On lit dans les \emph{Lois} 654a : 
 \og Les divinités dont nous avons parlé nous ont été données pour nous accompagner dans nos ch\oe urs de danse, et ce sont elles qui, à leur tour, nous ont donné le sentiment du rythme et de l'harmonie dans son union avec du plaisir\fg \cite[p. 674]{Lois}. 
 
 C'est surtout dans le \emph{Timée} que la musique est présentée comme élément divin lié à l'harmonie cosmique. Dans la section suivante, on développera ce thème.      
           
Le thème de l'harmonie des sphères qui est né avec les Pythagoriciens, fut développé par plusieurs auteurs de l'Antiquité. On peut citer à ce sujet Nicomaque,  le mathématicien grec du $I^{{\mathrm{er}}}$-II$^{{\mathrm{e}}}$  siècles, né à Gérase (en Palestine) qui,  dans son \emph{Manuel d'Harmonique} \cite{Nicomaque}, assigna aux corps célestes des notes musicales, une tradition qui fut reprise par le philosophe latin Boèce (V$^{{\mathrm{e}}}$-VI$^{{\mathrm{e}}}$  siècles) qui, inspiré par Nicomaque, fit de l'harmonie du monde, sous le nom de \emph{musica mundana}, une partie de son \emph{Institution musicale}\footnote{On doit à Boèce une division en trois parties de la musique :  \emph{musica mundana},  \emph{musica humana} et \emph{musica instrumentalis}, une tradition qui demeura vivante jusqu'à la Renaissance.}, un ouvrage qui demeura très influent durant tout le Moyen Âge européen. Ce thème fut poursuivi jusqu'après la Renaissance.   
  
  \section{Mathématiques et musique dans le \emph{Timée}}\label{s:Timee}

Parmi les écrits de Platon, le \emph{Timée} est l'un des derniers, et il y occupe une place particulière.  Il est considéré comme son ouvrage majeur, si ce n'est l' ouvrage philosophique majeur dans l'absolu. Il fut l'objet d'interprétations et de débats qui durent encore aujourd'hui, depuis deux millénaires. On ne compte plus le nombre de commentaires sur  le \emph{Timée}.  C'est l'ouvrage de Platon où les  relations de la philosophie avec les mondes des mathématiques et de la musique sont développées de façon explicite.   C'est en même temps un ouvrage de cosmologie, et, par conséquent, de physique, dans lequel Platon décrit la création de l'univers,  même s'il considérait qu'il ne peut y avoir réellement de science physique (ou science naturelle), puisque, disait-il, la nature est dans un état permanent de transformation, et qu'on ne peut décrire ce qui change.
   
      Le \emph{Timée} est un dialogue entre quatre personnages dont le sujet est la création du monde  (macrocosmos)  et celle de l'homme  (microcosmos). La création est présentée sous la forme d'une métaphore, racontée par Timée, celui parmi les quatre interlocuteurs qui représente les Pythagoriciens. Le récit n'est pas dogmatique, et Platon, comme il le dit lui-même, nous y mène dans le domaine de la fable, ou, tout au plus, celui du vraisemblable.     
L'un des thèmes principaux abordés est celui de  la musique et des harmonies, célestes, mais aussi humaines. Dans la cosmologie de Platon, le monde, tout comme l'homme, possède une âme et un corps. L'âme du monde est construite suivant des proportions, et ces proportions sont celles qu'un musicien d'aujourd'hui peut facilement reconnaître comme étant celles d'une gamme. Le corps du monde est construit à l'aide des éléments de base des polyèdres platoniciens et dans ces constructions, les proportions jouent un rôle central. Il faut rappeler ici que le mot grec pour proportion est \tg{l'ogos} (logos), que ce mot est ambigu car il a plusieurs significations, et Platon, comme d'autres philosophes grecs, comme pour entretenir l'ambiguïté autour de ce mot, l'utilise selon plusieurs sens en même temps. Le mot logos a des connotations mathématiques, logiques, éthiques, linguistiques, philosophiques et théologiques, et les grands auteurs grecs les ont toutes utilisées. Parmi les significations de ce mot, il faut mentionner, en plus de celui de  proportion, celui de fraction, de verbe (parole), de raisonnement, ou bien de raison (au sens d'une personne raisonnable).  Plus tard, le mot logos est entré dans la littérature chrétienne. Le début du \emph{Prologue} de Saint Jean le Théologien est scandé par ce mot :  \og Au commencement était le logos, et le logos était avec Dieu, et le logos était Dieu\fg. Les commentateurs n'ont pas manqué de souligner la référence au Démiurge de Platon, qui a créé le monde à partir du logos.
 Théon, dans son \emph{Exposition}, 
 consacre une section complète à ce mot, intitulée  \emph{En combien de sens se prend le mot logos}. Il y écrit:  \og On appelle [logos] le langage que les modernes nomment oral et le raisonnement mental sans émission de voix ; on appelle encore ainsi le rapport de proportion, et c'est en ce sens qu'on dit qu'il y a rapport de telle chose à telle autre; l'explication des éléments de l'univers; le compte des choses qui honorent ou qui sont honorées, et c'est dans cette acception que nous disons : tenir compte de quelque chose, ou n'en pas tenir compte. On appelle encore logos le calcul des banquiers, les discours de Démosthène et de Lysias dans leurs \oe uvres écrites; la définition des choses, qui en explique l'essence, puisque c'est à cela qu'elle sert; le syllogisme et l'induction; les récits libyques  et la fable. On donne aussi le nom de logos à l'éloge et au proverbe. C'est encore ainsi qu'on appelle la raison de la forme, la raison séminale et beaucoup d'autres.
Mais, selon Platon, on emploie le mot logos en quatre sens : \og On appelle ainsi la pensée mentale et sans parole, le discours procédant de l'esprit et exprimé par la voix, l'explication des éléments de l'univers et la raison de proportion.\fg  \cite{Theon} Ce passage est suivi d'une exposition des notions de fraction et de proportion  telles qu'elles sont utilisées en mathématiques, notamment celles qui interviennent en musique, et tout particulièrement celles qui sont utilisées dans le \emph{Timée}.

\begin{figure}[htbp]
\centering
\includegraphics[width=10cm]{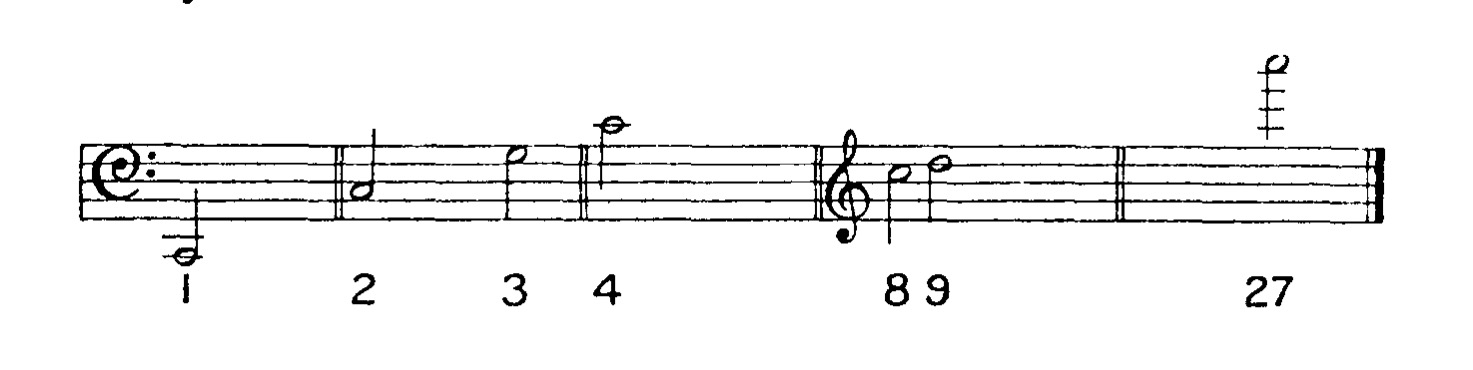}
\includegraphics[width=9.5cm]{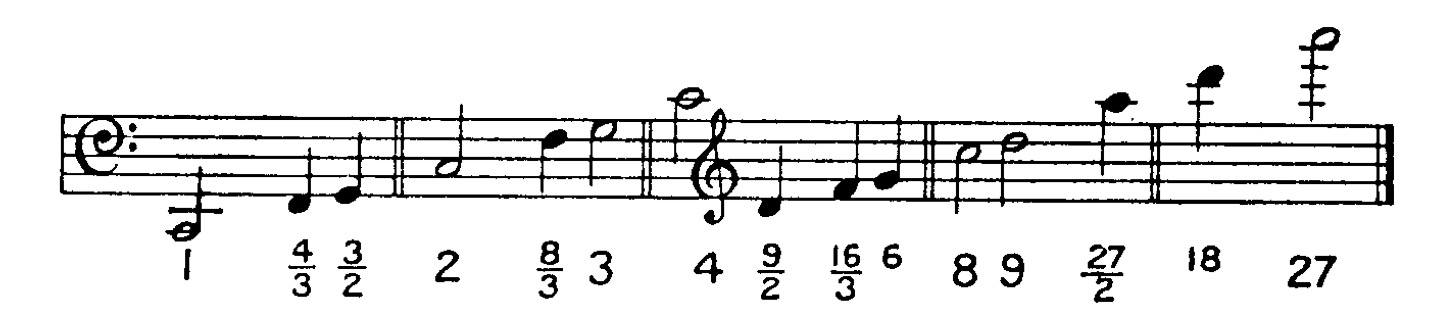}
\includegraphics[width=6cm]{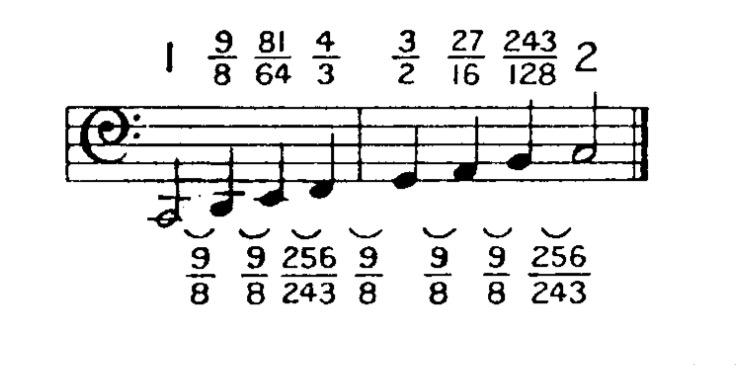}
\caption{\smaller } 
\label{gamme}
\end{figure}

Parmi les thèmes musicaux qui émergent de la construction du monde d'après Platon, le plus important est celui de la construction d'une gamme longue de quatre octaves et une sixte majeure, une gamme d'amplitude 27. Les proportions  (\tg{l'ogoi}, logoi) qui y interviennent sont celles suivant lesquelles l'âme du monde est construite à partir de trois substances : une substance qui est toujours identique à elle-même, désignée par \og Le Même\fg, une substance divisible désignée par \og L'Autre\fg, ou \og Le Différent\fg, et une troisième substance, qui est intermédiaire et qui sert à lier ensemble les deux premières. Il est hors de question de tenter d'expliquer ici la signification philosophique de ces constituants ; ce sont seulement les proportions qui nous importent dans ce contexte. Citons un extrait d'un passage parmi ceux qui ont servi à établir l'analogie musicale, un thème  que tous les commentateurs du \emph{Timée} ont souligné. Platon écrit : \og  Il se mit donc à faire les divisions que voici : il y eut en premier lieu une part qu'il préleva sur le tout ; après celle-ci, il en préleva une seconde, double de la première ;  puis une troisième, qui valait une fois et demie la seconde, était le triple de la première ; la quatrième fut double de la seconde, la cinquième triple de la troisième, la sixième égale à huit fois la première, et la septième à vingt-sept fois la première.   Après quoi, il combla les intervalles doubles ainsi que les triples du mélange, détachant encore des parts et les intercalant entre les premières, de sorte que dans chaque intervalle il y eût deux médiétés [\ldots] \fg (35b-36a, \cite[p. 450]{Timee}). Les médiétés dont il s'agit ici sont les moyennes  arithmétique et harmonique, celles qui sont utilisées pour remplir les intervalles musicaux. Cornford, dans son ouvrage \cite{Cornford}, décrit, en utilisant les notations de la musique moderne, la gamme construite à partir des nombres et des moyennes fournis par Platon. Il suit pour cela une description donnée par Archer-Hind dans une édition du Timée  publiée en 1888 \cite{AH}. Dans la figure \ref{gamme} nous avons reproduit un dessin décrivant trois étapes de cette construction, extrait de l'ouvrage de Cornford. La première ligne représente les sept notes avec lesquelles le Démiurge commence sa construction.  Les nombres qui leur sont associés représentent l'étendue de l'intervalle que la note fait avec la première de ces sept notes. Dans la deuxième ligne, les moyennes harmonique et arithmétique ont été introduites entre les termes consécutifs de cette série. Par exemple, $\displaystyle\frac{3}{2}$  est la moyenne arithmétique\footnote{La moyenne arithmétique de deux nombres $x_1, x_2$ est $\displaystyle\frac{1}{2}(x_1+x_2)$.} de $1$ et $2$ (car $\displaystyle\frac{3}{2}= \frac{1+2}{2}$) et $\displaystyle\frac{4}{3}$  est leur moyenne harmonique\footnote{La moyenne harmonique de deux nombres $x_1, x_2$ est $2\times\displaystyle\frac{x_1\times x_2}{x_1+x_2}$.} (car $\displaystyle\frac{4}{2}= 2\frac{1\times 2}{1+ 2}$). La troisième ligne représente le résultat obtenu en remplissant la première octave qui apparaît dans cette gamme, présentée comme une juxtaposition de deux tétracordes,  par des intervalles de ton (valeur $\displaystyle\frac{9}{8}$) et un petit intervalle dont la valeur ($\displaystyle\frac{256}{243}$) est proche de celle d'un demi-ton.

      Les références à la musique, dans la construction du monde, vont au-delà de ce qui se rapporte à la gamme. Il y a divers passages qui évoquent le rythme, le chant, la voix et d'autres notions liées au son et à la musique. Par exemple, dans le \emph{Timée}  (47c, \cite[p. 465]{Timee}), Platon écrit : \og Pour ce qui est de la voix et de l'ouïe, derechef il en faut dire autant : c'est aux mêmes fins et pour les mêmes motifs que les dieux nous en ont fait présent. La parole en effet a été disposée précisément pour ces fins et y apporte la plus large contribution ; et  dans la mesure où la voix est utilisée en musique et s'adresse à l'ouïe, c'est en vue de l'harmonie qu'elle nous est donnée. Or, l'harmonie est faite de mouvements de même nature que les révolutions de l'âme en nous ; et pour qui use avec l'intelligence du commerce des Muses, ce n'est pas dans un plaisir irraisonné, comme le veut l'opinion actuelle, que réside son utilité ; mais comme, de naissance, en nous la révolution de l'âme est inharmonique, c'est pour la mettre en ordre et en accord avec soi que l'harmonie nous a été donnée pour alliée par les muses. Et le rythme, à son tour, c'est à cause d'une absence en nous de mesure, et d'un manque de grâce que chez la plupart manifeste le maintien, que les mêmes divinités nous l'ont à cette fin donné comme remède.\fg

\begin{figure}[htbp]
\centering
\includegraphics[width=8.51cm]{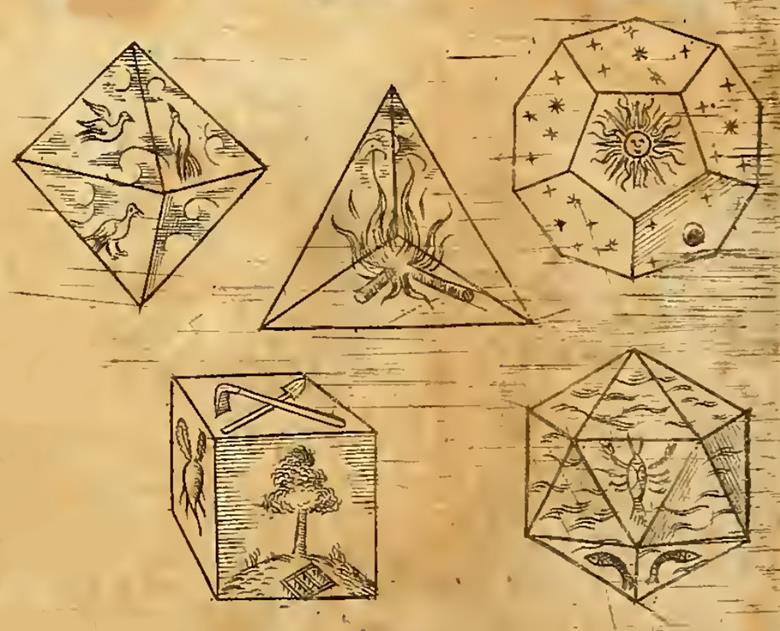}
\includegraphics[width=3.35cm]{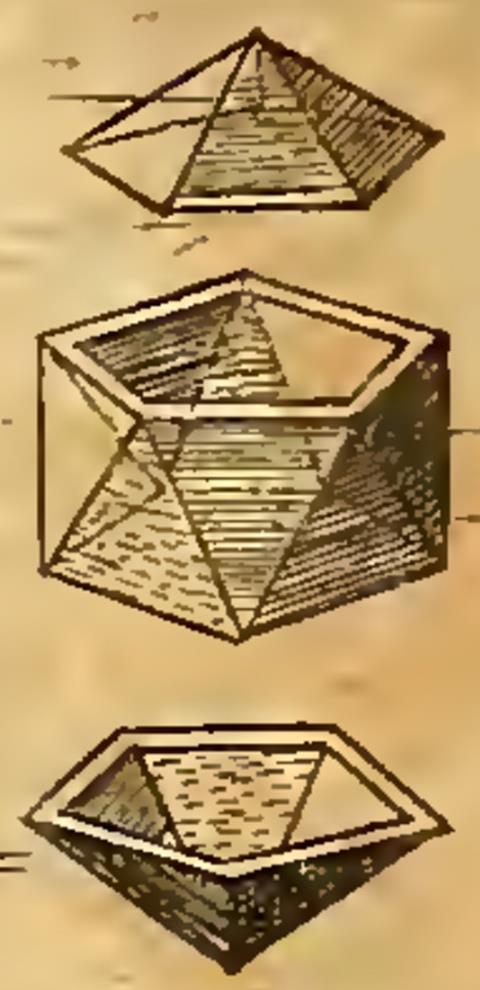}
\caption{\smaller Extrait du \emph{Harmonices Mundi} de Kepler (1619), représentant les cinq polyèdres platoniciens. Dans la figure de droite, Kepler a dessiné une décomposition de l'icosaèdre.} 
\label{Kepler1}
\end{figure}

Parmi les aspects géométriques du \emph{Timée}, il y a celui de la description des polyèdres platoniciens, dont je voudrais parler encore ici. Dans le passage du \emph{Timée} où le Démiurge construit le corps du monde, Platon s'étend sur ces polyèdres et leur construction. Ce corps est bâti en utilisant les quatre éléments premiers de la matière, ceux qu'Empédocle et d'autres Présocratiques ont dégagé :  la terre, l'eau, l'air et le feu\footnote{Il est intéressant de noter que le mot grec  \tg{stoiqe\~ia} (sticheia) que les Présocratiques ont utilisé pour désigner les éléments premiers de la matière est le même que celui qu'Euclide utilisa pour les éléments de la géométrie.}. Platon associe ces éléments indivisibles (\og atomes\fg) au cube,  icosaèdre, octaèdre et tétraèdre respectivement. On a beaucoup écrit sur les raisons de cette association.  
On peut voir ces polyèdres dans la figure \ref{Kepler1} ci-dessous, tirée de l'\emph{Harmonices Mundi} de Johannes Kepler,  l'un des principaux fondateurs de la science moderne, qui se déclarait lui-même héritier de Platon.
Les commentateurs de Platon, comme Kepler, ont généralement considéré que le cube est associé à la terre parce que, posé sur un plan, il est plus stable que les autres polyèdres ; que l'icosaèdre est parmi les polyèdres réguliers celui qui a le plus grand nombre de faces, et que par conséquent sa forme se rapproche, plus que celle des autres, d'une forme ronde qui évoque des gouttes de liquide ; que 
l'octaèdre, quand il est placé dans la position qu'il a dans la figure \ref{Kepler1}, indiquant quatre directions au niveau supérieur et quatre au niveau inférieur, évoque l'air ; et que le tétraèdre, comme objet pointu et ayant un petit nombre de faces,  évoque le feu\footnote{Une analyse philosophique originale et intéressante des passages du \emph{Timée} évoquant les éléments de la matière est donnée par Negrepontis et Kalisperi dans leur article récent \cite{NK}. Ces auteurs fournissent en particulier une interprétation complètement nouvelle, basée sur une lecture attentive, du passage 48b5-c2 du \emph{Timée}. Ils expliquent que les objets sensibles, chez Platon, ne peuvent pas être des \og principes\fg ou bien des \og éléments\fg (\tg{stoiqe\~ia}), mais qu'il faut les considérer comme des \emph{syllabes} (des dyades, triades, etc.), la base de cette perpective découlant de l'interprétation nouvelle et mathématique du langage platonicien donnée par Negrepontis dans \cite{Negrepontis} suivant laquelle un objet sensible participe à un \^Etre intelligible comme une dyade ayant une anthyphairesis égale à un sous-segment initial fini de l'anthyphairesis de l'\^Etre intelligible. En particulier, un \^Etre sensible n'est pas un objet unique, comme une surface ou un polyèdre régulier, mais une paire d'objets. Signalons aussi l'article \cite{N2} de Negrepontis où ce dernier donne une nouvelle interprétation du langage utilisé par Platon dans la section 142b-155e du \emph{Parménide}, en termes d'un analogue philosophique de l'anthyphairesis dans lequel la théorie de la musique  joue un rôle fondamental.}.

 Dans un souci de faire en sorte que tout dérive d'un certain nombre d'éléments premiers, Platon explique que les faces du cube sont constituées de quatre triangles isocèles rectangles, et que celles du tétraèdre sont des triangles équilatéraux, qui, divisés par une hauteur, donnent pour chacun deux triangles rectangles dont l'hypoténuse est le double du côté le plus petit. Des divisions des faces des autres polyèdres sont décrites de la même manière. Platon dit que même s'il y  a une infinité de formes possibles pour les triangles, seulement deux d'entre elles sont nécessaires pour construire le corps du monde : le triangle rectangle isocèle, et celui qu'on obtient en coupant suivant une hauteur un triangle équilatéral. Ces deux triangles jouent le même rôle que les ``éléments" de la matière. Le plus beau parmi les triangles, dit Platon,  est le triangle équilatéral (54a, \cite[p. 474]{Timee}). Il écrit : \og Des deux triangles dont nous parlions, celui qui est isocèle n'a qu'une seule nature en partage ; celui qui est allongé en a une infinité ; il nous faut donc, dans cette infinité encore, préférer le plus beau, si nous entendons commencer suivant le procédé voulu. Si donc quelqu'un en peut, d'après son choix, indiquer un plus beau  pour la constitution de ces corps, nous ne verrons pas en lui un adversaire, mais un ami, et reconnaîtrons sa victoire. Nous posons que, de cette multitude de triangles, il en est un qui est le plus beau de tous, et nous passons par-dessus tous les autres : c'est celui dont deux en constituent un troisième, le triangle équilatéral. \fg  Une discussion suit, mais ce que je voulais souligner ici c'est l'argument esthétique, et là aussi l'approche de Platon est celle d'un mathématicien, pour qui le critère final est celui de la beauté du résultat. 
 
Les analogies se poursuivent dans le le passage qui commence en 55a, dans lequel Platon  fait un décompte du nombre de triangles rectangles isocèles ou équilatéraux nécessaires pour construire chacun des polyèdres.  
\`A partir des propriétés des faces des polyèdres qui sont convertibles l'une dans l'autre, Platon déduit quels sont les éléments de la matière qui sont transformables l'un dans l'autre  (\emph{Timée} 53e-54b) : air, feu, eau peuvent se transformer l'un dans l'autre, car ils ont le même type de faces, etc. Le cinquième polyèdre régulier est aussi mentionné dans le \emph{Timée} (55c), même s'il n'est associé à aucun des quatre éléments : Platon dit que le Démiurge l'utilisa pour \og le tout\fg, un mot qui est généralement interprété comme désignant la voûte céleste, qui a la forme d'une sphère\footnote{\label{f:Plutarque}Plutarque, dans ses \emph{Questions platoniques} (1003d), note que parmi les cinq polyèdres  platoniciens, c'est le dodécaèdre qui s'approche le plus de la sphère. Il écrit : \og [\ldots] Car la multitude des éléments du dodécaèdre et la grande ouverture de ses angles font que, s'éloignant beaucoup de la ligne droite, il se courbe facilement, et son périmètre, comme dans les sphères composées de douze pièces réunies, approche davantage de la forme circulaire et contient un très grand espace. Il y a vingt angles solides, dont chacun est renfermé dans trois angles plans et obtus qui contiennent chacun un angle droit et la cinquième partie de cet angle. D'ailleurs le dodécaèdre est formé de douze pentagones, dont les côtés et les angles sont égaux, et composés chacun des trente premiers triangles scalènes. Il semble donc être une image du zodiaque et de l'année, puisque ses divisions sont égales à l'un et à l'autre.\fg}, le nombre de ses faces étant le même que celui des douze signes du zodiaque. Kepler n'a pas manqué de souligner cette comparaison (voir la figure \ref{Kepler1}). Dans cette     
   figure, qui représente les cinq solides platoniciens et qui est tirée de l'\emph{Harmonices Mundi} de Kepler, quatre de ces polyèdres sont associés aux éléments de la nature, comme Platon l'a con\c cu,  et le cinquième, le dodécaèdre, est associé aux douze signes du zodiaque (comme l'avait pensé Plutarque avant lui, voir la note \ref{f:Plutarque}).     
  Proclus et les autres commentateurs d'Euclide et de Platon appelaient les polyèdres platoniciens les \emph{figures cosmiques}, plusieurs siècles avant que Kepler en fit la base de ses réflexions géométriques et astronomiques. On en reparlera plus loin.
 
    La figure du Lambda (Figure \ref{Macrobe}) réunit l'arithmétique et la géométrie du triangle.
    
  Je voudrais  terminer cette section sur les éléments de la matière en citant un texte du physicien Werner Heisenberg. Il s'agit d'un extrait de l'ouvrage  \emph{Der Teil und das Ganze: Gespräche im Umkreis der Atomphysik} (La partie et le tout: conversations dans le domaine de la physique atomique) \cite{H2} dans lequel l'auteur raconte comment, durant l'été 1919,  il tomba sur des passages du \emph{Timée} qui l'avaient déjà fasciné auparavant, précisément les passages dans lesquels Platon parle des particules élémentaires de la matière, disant que ces particules élémentaires sont  des triangles rectangles pouvant être assemblés pour obtenir des triangles isocèles et des carrés, qui, ensuite, se combinent en des polyèdres réguliers. 
Cet extrait du livre de Heisenberg vient, dans l'ouvrage, avant que l'auteur n'entame la discussion sur le problème de l'atome\footnote{Cf. p. \cite{H2} et ss. de la traduction anglaise.}. 

Heisenberg déclare que la lecture du passage de Platon sur les parties constituantes géométriques de la matière fut déterminante pour lui. Il écrit que sa première réaction, après cette lecture, fut de ne pas pouvoir se décider sur la question de savoir si Platon voulait dire que ces polyèdres étaient associés aux quatre éléments de la matière (terre, feu, air et eau) uniquement comme symboles, ou si plutôt, il entendait par là que les plus petites parties de l'élément terre avaient des formes cubiques, celles du feu tétraédriques, etc. Bien qu'au premier abord, dit-il, il lut ce texte avec une certaine suspicion, et qu'il fut perturbé par le fait qu'un grand philosophe comme Platon affirmât de telles fantaisies, cela lui donna l'impulsion pour réfléchir à la question, et pour rechercher un principe qui pourrait justifier ces spéculations. Il fut de plus en plus captivé, dit-il, par l'idée que les plus petites particules de matière se réduisent à des formes géométriques, et par le fait que Platon ait choisi comme modèles géométriques les polyèdres réguliers. Il  se posa la question de savoir pourquoi Platon pouvait voir de l'ordre dans la nature, alors que lui ne le pouvait pas, et il en arriva ensuite à se poser la question trouver de toutes les significations possible du mot ordre, et celle de trouver un point central qui pouvait réunir toutes les notions d'ordre fragmentaire. Heisenberg raconte que quelques jours plus tard, il entendit un violoniste jouer la Chaconne de Bach. Citons les mots de Heisenberg :  \og Les phrases claires de la chaconne,\fg écrit-il, \og me touchèrent comme un vent froid, brisant la brume et révélant les structures imposantes au-delà. Il y a toujours eu dans le langage de la musique, dans la philosophie et dans la religion, un chemin vers l'ordre central, aujourd'hui pas moins que du temps de Platon et de Bach. Je savais cela maintenant de par ma propre expérience. \fg 

Le nom de Heisenberg est attaché au principe d'indétermination en mécanique quantique, un principe qui affirme qu'il est impossible de déterminer avec précision deux valeurs simultanées de deux propriétés physiques d'une même particule. Un tel principe peut finalement être mis en parallèle avec l'affirmation de Platon suivant laquelle on ne peut pas connaître ce qui change en permanence.

       \section{Cercles}\label{s:cercles}

  Je voudrais maintenant parler d'une figure géométrique qui intervient de manière récurrente en musique, depuis l'Antiquité jusqu'à nos jours ; il s'agit du cercle. Considéré par Platon comme l'une des figures parfaites, à cause de ses  symétries, le cercle est le lieu du mouvement des astres. De plus, c'est l'analogue uni-dimensionnel de la sphère, celle qui, d'après le \emph{Timée}, englobe l'univers tout entier : dans le \emph{Timée} 33b, Platon écrit que le Démiurge donna à l'univers la forme ronde d'une sphère, une figure qui, dit-il, est, comme le cercle, semblable à elle-même en tous ses points, ce qui en fait la figure de l'espace \og la plus juste qui soit\fg. Cette justice est aussi symbolisée par le fait que tous ses points sont équidistants d'un point donné (le centre), et par le fait qu'ils sont tous semblables : contrairement aux points d'un solide platonicien, il n'y a, sur la sphère, ni sommet, ni arête, ni face. Utilisant un langage mathématique, cette propriété dit qu'on peut passer de n'importe quel point de la sphère à n'importe quel autre par une rotation.
  Dans le passage 
     34a-b (\cite[p. 449]{Timee}), on lit  à propos du corps du monde  : \og  Tout calculé, [le Démiurge] le fit bien poli, sans inégalités dans sa surface, en tous ses points équidistant du centre ; ce fut un tout, un corps complet, fait de corps au complet.\fg Dans le même passage, à propos de l'âme du monde, on lit :  \og Il la plaça au centre du monde, puis l'étendit à travers toutes ses parties et même en dehors, de sorte que le corps en fut enveloppé ; cercle entraîné dans une rotation circulaire, c'est là comme il établit le Ciel : rien qu'un seul, solitaire, capable en vertu de son excellence d'être en union de soi à soi sans avoir besoin de rien d'autre, objet de connaissance et d'amitié pour soi-même, à en être comblé ! C'est par tous ces moyens qu'il le fit naître  Dieu bienheureux.
 \fg

Les cercles et les rotations circulaires apparaissent dans plusieurs autres passages du  
 \emph{Timée}, ce qui n'est pas étonnant, vu le contenu astronomique du traité. Par exemple, dans le passage 36b-c \cite[p. 451]{Timee}, dans lequel Platon décrit la construction de l'âme, il déclare qu'après avoir ordonné les mélanges faits suivant des proportions (musicales) sur deux lignes droites, le Démiurge plaça ces deux lignes l'une sur l'autre, sous la forme d'une lettre X, (la lettre grecque $\chi$) puis il les courba en cercles, réunit les deux extrémités de chacun des segments, et il leur imprima le mouvement perpétuel du cercle. Il assigna à la nature du Même le mouvement du cercle extérieur, et à celle de l'Autre celui du cercle intérieur.  De nouveau, les cercles sont divisés suivant des proportions qui, on le sait maintenant,  sont celles qui sont utilisées pour construire les gammes musicales :  \og [\ldots] Quant à la révolution intérieure, il la divisa en six pour en faire sept cercles inégaux, dont les distances correspondent aux intervalles doubles et triples, ceux-ci étant trois de chaque série : il prescrivit à ces cercles d'aller en sens inverse les uns des autres, trois avec des vitesses égales, les quatre autres avec des vitesses inégales entre elles et à celles des trois précédents, mais tous d'un mouvement bien réglé.\fg
 Dans le passage 44a-44b \cite[p. 461]{Timee}, Platon parle des \og cercles de l'âme\fg : l'âme, comme le monde, a deux cercles : celui du Même, qui est son équateur, et celui de l'Autre, qui est son écliptique.  Dans le passage 47b-c \cite[p. 465]{Timee}, il écrit: 
   \og  [\ldots] Si un dieu a pour nous inventé le présent de la vue, c'est afin que, contemplant au ciel les révolutions de l'intelligence, nous en fassions application aux cercles que parcourent en nous les opérations de la pensée ; ceux-ci sont de même nature que celles-là. \fg

\begin{figure}[htbp]
\centering
\includegraphics[width=6cm]{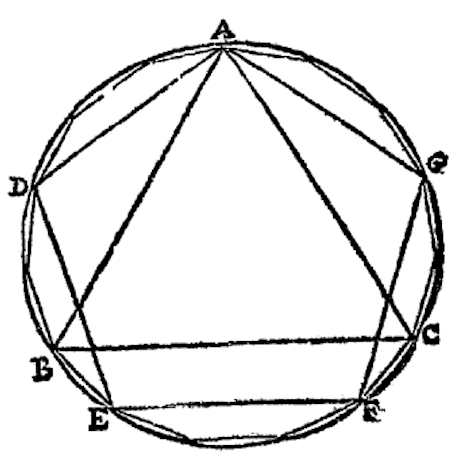}
\includegraphics[width=6cm]{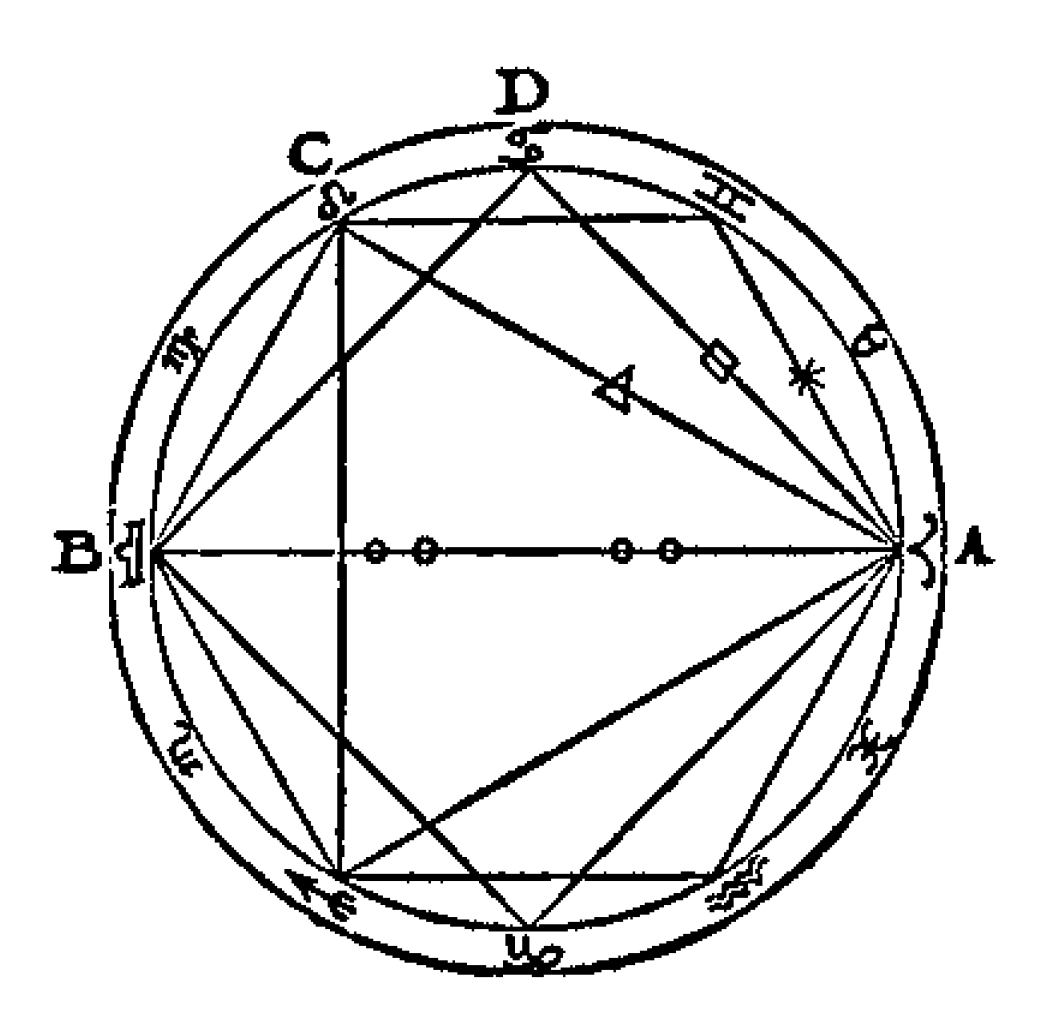}
\includegraphics[width=6cm]{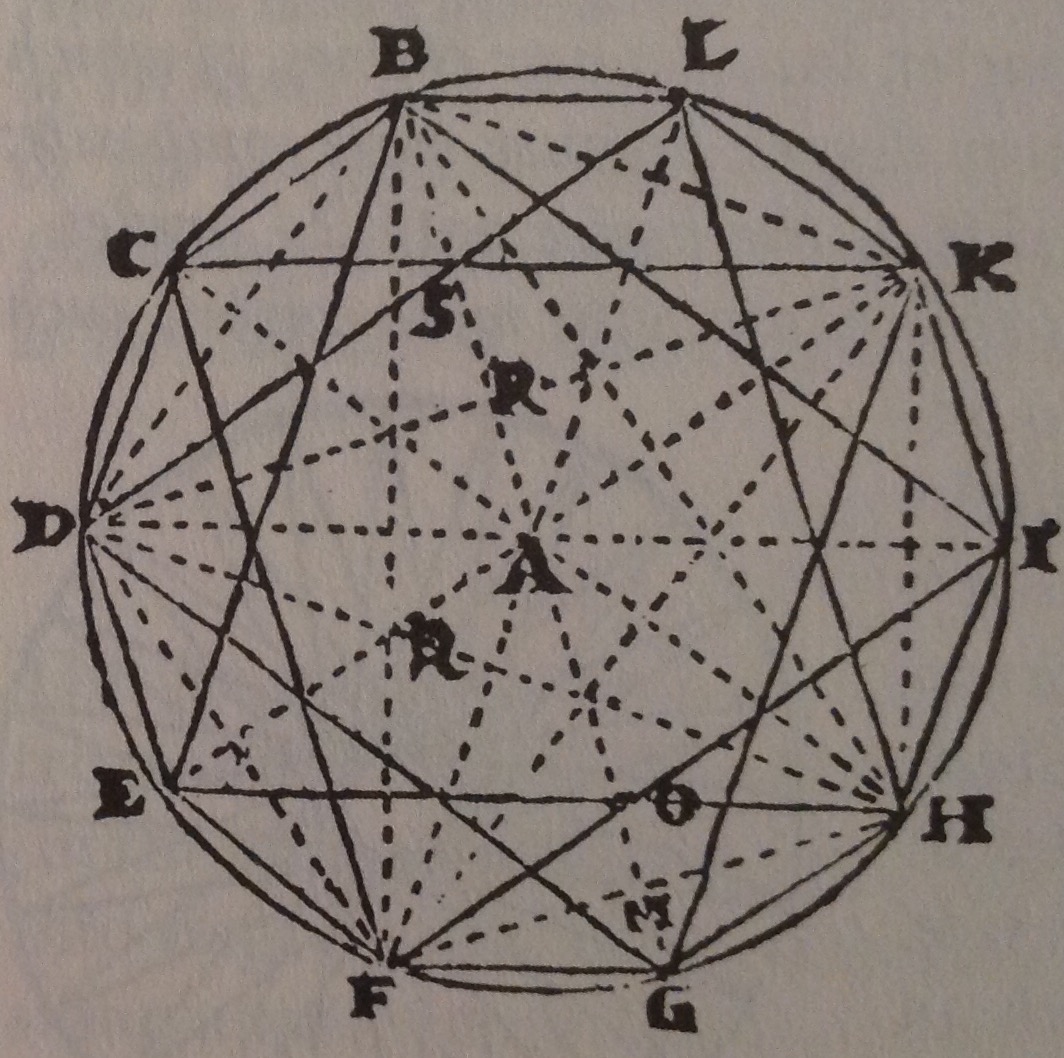}
\includegraphics[width=6.2cm]{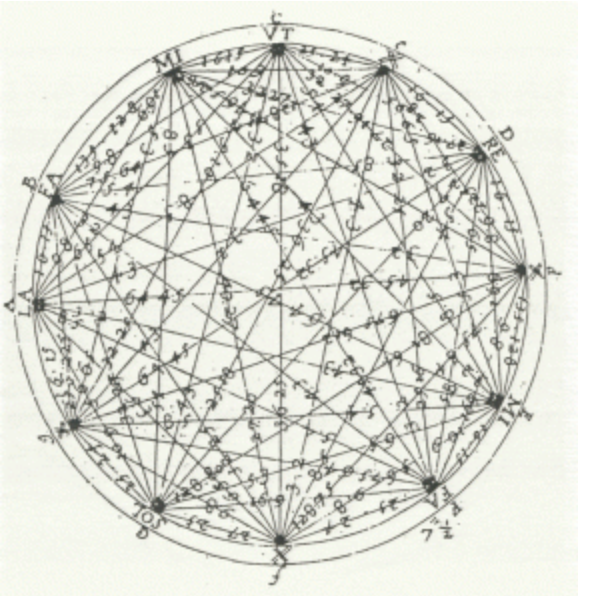}
\includegraphics[width=6cm]{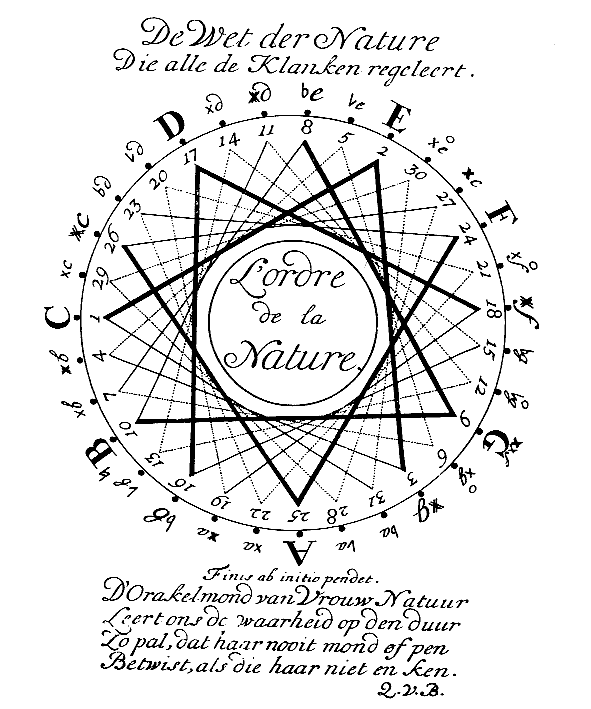}
\includegraphics[width=5.5cm]{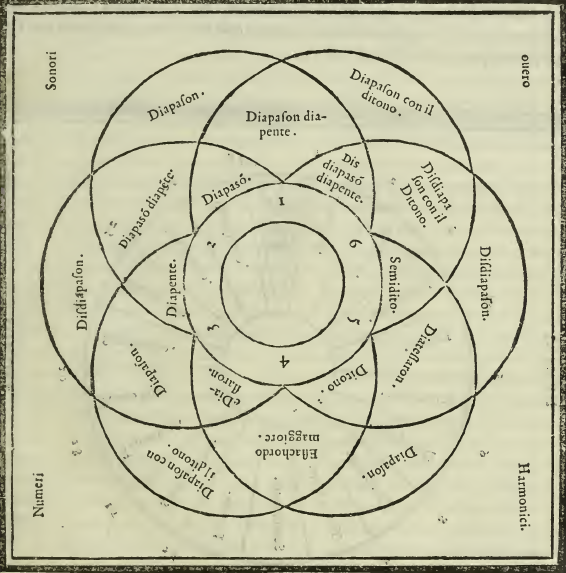}
\caption{\smaller De gauche à droite et de haut en bas: Extrait des \emph{\'Eléments} d'Euclide (Proposition 16 livre IV) (vers 300 av. J.-C.),  de l'\emph{Harmonique} de Ptolémée (2$^{e}$ s. ap. J.-C.), de l'\emph{Harmonices Mundi} de Kepler (1619), de l'\emph{Harmonie universelle} de Mersenne (1636), des \emph{Institutioni harmoniche} (1558) de Zarlino et des  \emph{Elementa musica} de Quirinus van Blankenburg (1739)} 
\label{f:cercles-Ptolemee}
\end{figure}

Dans la figure  \ref{f:cercles-Ptolemee}, nous avons reproduit des cercles, complétés par d'autres lignes, apparaissant en musique, en mathématiques et en astronomie, les trois sujets du \emph{Timée}. Les figures sont tirées des écrits d'Euclide, Ptolémée, Kepler, Mersenne,  Zarlino et van Blankenburg. Le dessin  tiré des \emph{\'Eléments} d'Euclide concerne la construction de certains polygones réguliers. Kepler développa longuement, dans le \emph{Harmonices Mundi}, la relation entre la constructibilité de ces polygones à la règle et au compas et les problèmes de consonance musicale. 
Ptolémée et Kepler sont parmi les principaux mathématiciens-astronomes de tous les temps qui développèrent une analogie entre la structure du système planétaire et les consonances musicales, ces deux phénomènes étant pour eux l'expression de la même harmonie de l'univers. Le cercle joue un rôle primordial dans leurs écrits sur ce sujet.

Dans le  chapitre II du livre III de l'\emph{Harmonices Mundi}, intitulé \og Sur les divisions harmoniques d'une corde\fg, Kepler déclare que dans le chapitre XII du \emph{Mysterium Cosmographicum}, l'ouvrage qu'il avait écrit vingt-deux ans plus tôt, il avait pensé, par erreur, que les causes des divisions et des harmonies se trouvaient dans les cinq solides platoniciens, mais que maintenant, il pense que les solides platoniciens et les harmonies musicales trouvent tous les deux leur origine dans les divisions du cercle par les polygones réguliers.
Dans l'introduction du Livre I de l'\emph{Harmonices Mundi}, il écrit : \og Nous devons rechercher les causes des rapports harmoniques dans les divisions du cercle en quelques parties égales qui sont faites géométriquement et de manière à ce qu'elles soient connues, c'est-à-dire à partir des figures planes régulières constructibles.\footnote{Pour ce passage et certains autres, j'ai légèrement remanié la traduction de \cite{Kepler-HM}.} \fg \`A plusieurs endroits de son \oe uvre, Kepler considère le cercle comme une corde courbe qui peut émettre des sons musicaux.  Le chapitre I du livre III du même ouvrage, intitulé \og Des causes des consonances\fg, commence par une série de définitions et d'axiomes. \`A la fin du commentaire sur l'Axiome I, on lit : \og   Pour ce qui concerne la musique, il suffit qu'une corde tendue comme une ligne droite puisse être divisée de la même manière que quand elle est courbée sous la forme d'un cercle et divisée par les côtés d'une figure inscrite.\fg 
Dans l'Appendice au Livre V du même ouvrage, Kepler écrit : 
\og Le cercle coupé géométriquement par l'inscription des figures planes définit véritablement et de façon appropriée la proportion harmonique par la comparaison des parties avec le tout, le cercle se trouve donc dans les âmes par une certaine raison essentielle, formellement et abstraitement, non seulement issue de la matière, mais encore issue d'une certaine quantité elle-même, considérée matériellement ; c'est pourquoi les harmonies se trouvent avec le cercle dans les âmes, et il y a cette cause pourquoi les âmes sont émues par les harmonies.\fg

Dans le même appendice, Kepler rappelle comment Ptolémée établit une comparaison entre le zodiaque et le système ou l'échelle musicale, au chapitre VIII du livre III de ses \emph{Harmoniques}, et une autre comparaison entre les proportions tirées des consonances avec celles qui sont liées à la disposition des planètes, au chapitre IX du même ouvrage (on renvoie au dessin tiré des \emph{Harmoniques} de Ptolémée reproduit ici dans la Figure \ref{f:cercles-Ptolemee}).

\begin{figure}[htbp]
\centering
\includegraphics[width=12cm]{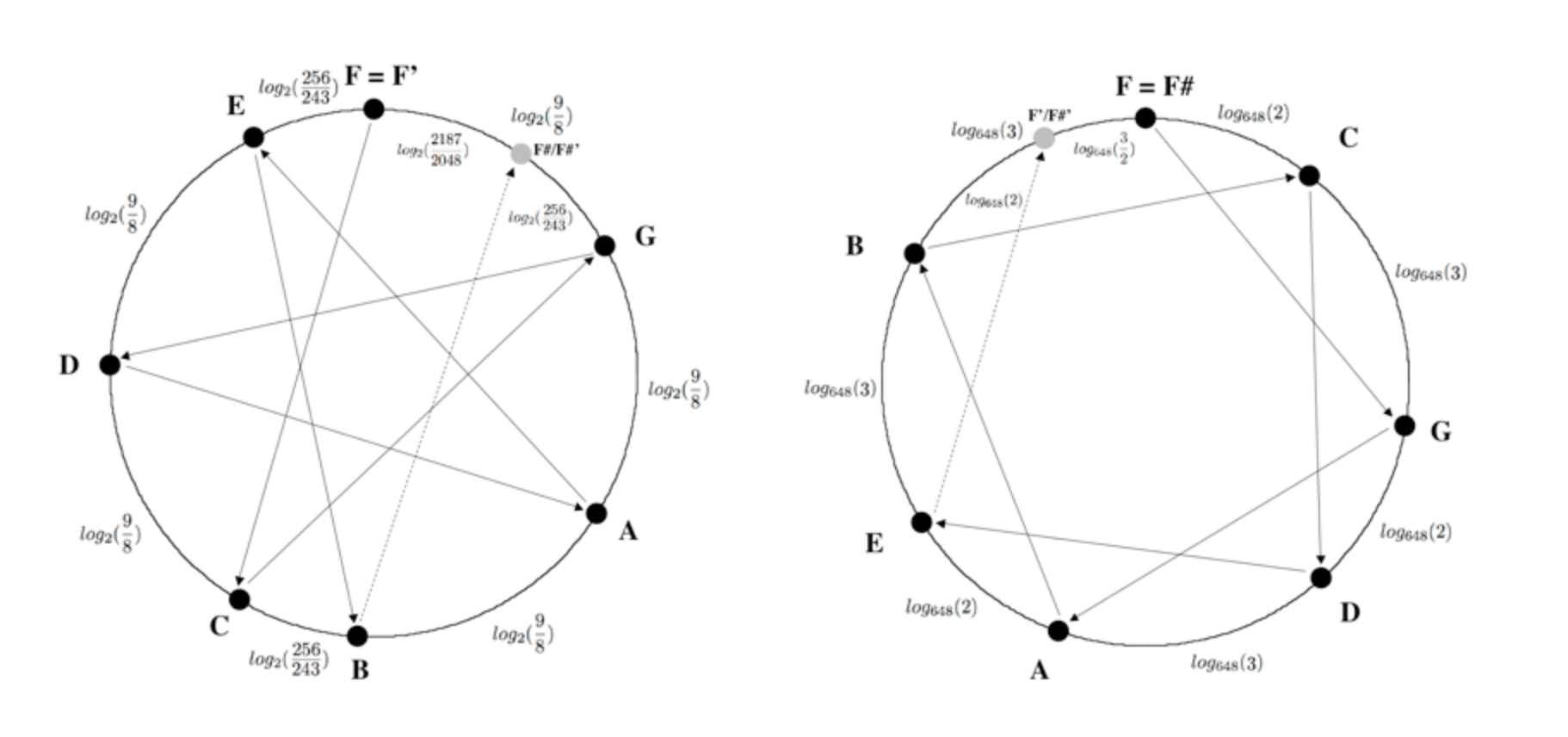}
\caption{\smaller Extrait d'un article de Thomas Noll.} 
\label{Noll}
\end{figure}

Cette tradition de représenter les intervalles musicaux sur des cercles se poursuit aujourd'hui. 
La Figure \ref{Noll} est tirée  d'un article récent de Thomas Noll \cite{Noll} ; il y a une multitude d'autres exemples que l'on pourrait citer.

    \section{Conclusion}
 
 Platon, comme tous les grands géomètres, était un poète et un rêveur. Il nourrissait son auditoire d'idées qui, à son époque, étaient complètement nouvelles. Il plaçait la musique et la géométrie à des niveaux de perfection et d'abstraction qui ne furent peut-être jamais dépassés après lui. Même s'il suggérait des analogies entre certains concepts mathématiques et le monde réel, ces analogies demeuraient théoriques, faisant de son discours, comme celui de Platon, une métaphore dans laquelle les mathématiques et la musique jouent les deux rôles principaux. Son \oe uvre n'a jamais cessé d'influencer la philosophie et la pensée occidentales.
 
  Le physicien Werner Heisenberg, que nous avons déjà mentionné, a exprimé à plusieurs occasions son admiration pour Platon, et en particulier pour le \emph{Timée}. Il déclara à plusieurs reprises que c'est la lecture de cet ouvrage qui lui donna l'idée d'entreprendre des études de physique. Nous avons déjà cité à la fin de la section \ref{s:Timee} un extrait de son livre \cite{H2} dans lequel il décrit comment la lecture du passage de Platon sur les particules élémentaires l'a influencé dans ses propres recherches. Le même ouvrage contient plusieurs autres références à Platon. On peut mentionner aussi l'article \emph{Platons Vorstellungen von den kleinsten Bausteinen der Materie und die Elementarteilchen der modernen Physik} (Les idées de Platon sur les plus petits éléments constitutifs de la matière et les particules élémentaires de la physique moderne)  \cite{H1} de Heisenberg.

Concernant l'impact de Platon sur la philosophie moderne, il faut mentionner les écrits d'Alfred North Whitehead, dont on cite souvent une phrase qui dit que la philosophie occidentale consiste en des notes de bas de page aux \'ecrits de Platon. La phrase est tirée de son ouvrage \emph{Process and Reality} \cite[p. 39]{Whitehead}, et je voudrais terminer cet article en citant le passage suivant qui contient cette phrase :

  \begin{quote}\small
  La  caractérisation générale la plus sûre de la  tradition philosophique européenne est  de dire qu'elle consiste en des notes de bas de page à Platon.
   Par là, je n'entends pas le schéma de pensée systématique que les érudits ont cru extraire de ses écrits. Je fais allusion plutôt à la richesse des idées générales qui y sont disséminées. Ses talents personnels, les possibilités très larges que lui offrait son expérience dans une grande époque de la civilisation, l'héritage qu'il a reçu d'une tradition intellectuelle qui n'a pas été encore par été raidie par une systématisation excessive, ont fait de ses idées une mine inépuisable d'inspiration\footnote{The safest general characterization of the European philosophical tradition is that it consists in a series of
footnotes to Plato. I do not mean the systematic scheme of thought which scholars have doubtfully extracted from his writings. I allude to the wealth of general ideas scattered through them. His personal endowments, his wide opportunities for experience at a great period of civilization, his inheritance of an intellectual tradition not yet stiffened by excessive systematization, have made his writing an inexhaustible mine of suggestion.}.
\end{quote} 
  
\begin{figure}[htbp]
\centering
\includegraphics[width=12cm]{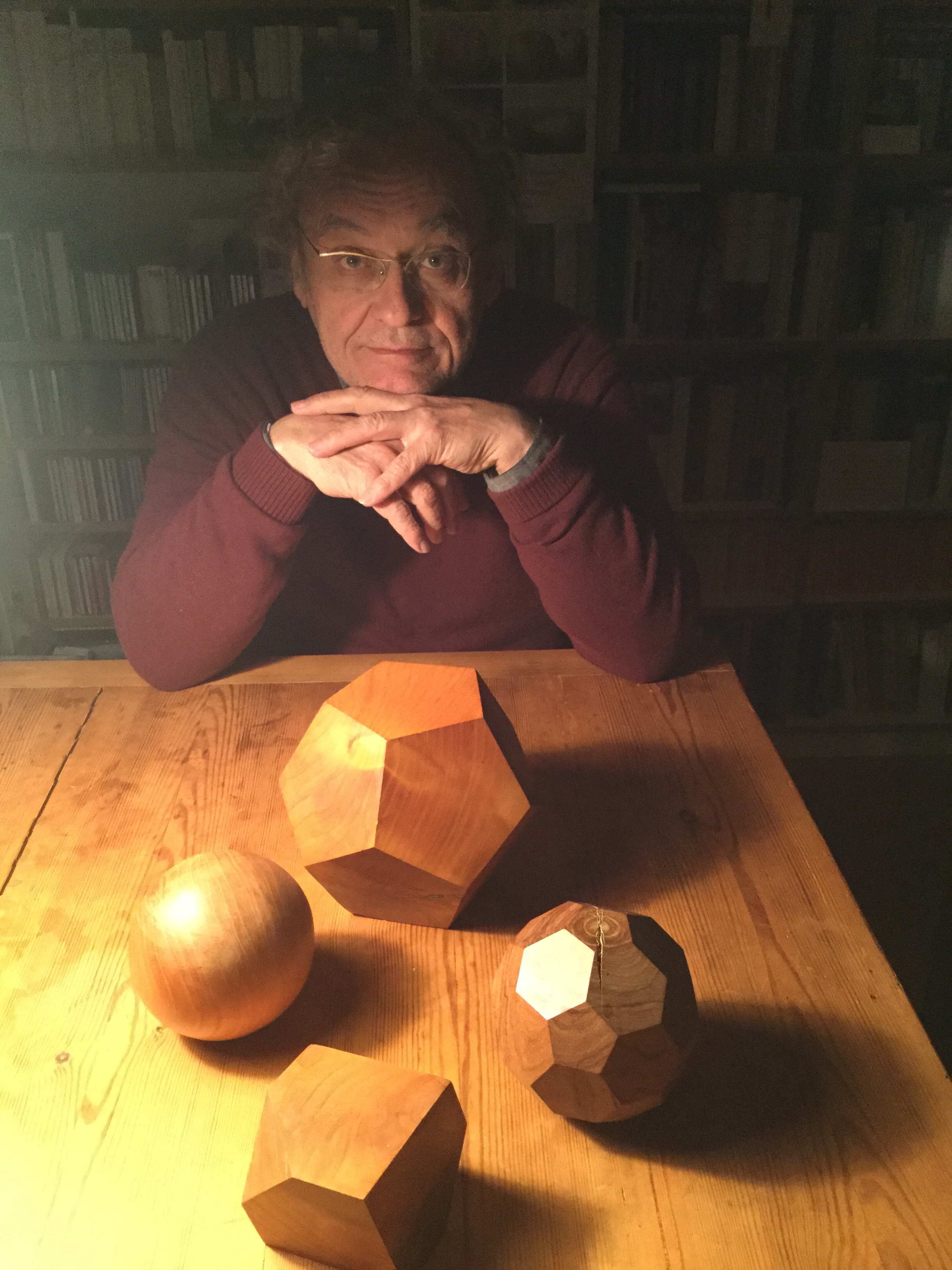}
\caption{\smaller L'apprenti essaye de trouver la position juste du polyèdre de D\"urer.
(Sculptures en bois de cerisier réalisées par Martin Turner.)} 
\label{Z1}
\end{figure}

\begin{figure}[htbp]
\centering
\includegraphics[width=12cm]{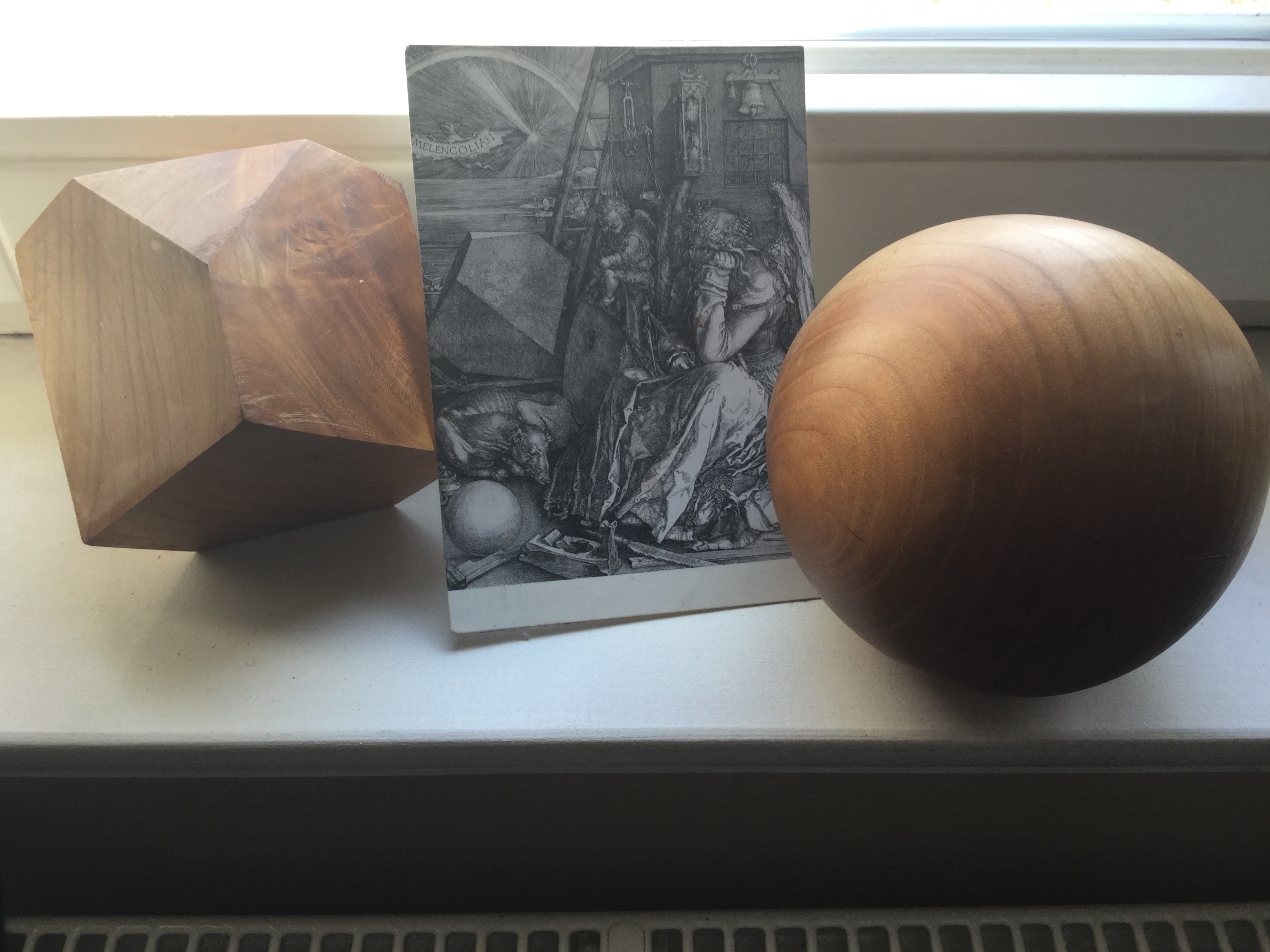}
\caption{\smaller Le problème résolu.} 
\label{Z2}
\end{figure}

  Et pour finir, deux photos : dans la première (Figure \ref{Z1}), on voit
Walter Zimmermann, en apprenti Platonicien---comme il se présente
lui-même---dans sa maison de campagne, devant une sphère et trois
polyèdres : le dodécaèdre régulier ; un polyèdre semi-régulier dont les
faces sont des pentagones et des hexagones réguliers ; et un polyèdre
non régulier à huit faces, dont les faces sont de deux types différents
(un triangle équilatéral et un pentagone non régulier). Ce polyèdre à
8 faces est celui que Dürer a représenté dans sa gravure Melancholia
(Figure \ref{Z2}).
  
   \bigskip

 \smaller \noindent {\bf Remerciements.} Je voudrais remercier Marie Corti et Marie-Pascale Hautefeuille pour les corrections qu'elles ont faites sur un première version de cet article.
\larger

  \vfill\eject


\begin{thebibliography}{999}
 
 \bibitem{AH} R. D. Archer-Hind (dir.) The Timaeus of Plato, Macmillan, London, 1988.
 
%
%
%




\bibitem{A-Belis} A. Belis, \emph{Aristoxène de Tarente et Aristote, Le traité d'harmonique}, coll. \'Etudes et commentaires, no. 110, Paris, Klincksieck,  1986.



\bibitem{Cornford} F. M. Cornford, Plato's cosmology: The \emph{Timaeus} of Plato, Hackett, Cambridge, 1935.



\bibitem{H1} W. Heisenberg, Platons Vorstellungen von den kleinsten Bausteinen der Materie und die Elementarteilchen der modernen Physik, In : Umkreis der Kunst. Eine Festschrift f\"ur E. Preetorius, Wiesbaden, 1953, p. 137-140.

\bibitem{H2} W. Heisenberg,  Der Teil und das Ganze: Gespräche im Umkreis der Atomphysik, Munich, 1969. Trad. anglaise soue le titre: Physics and beyond, encounters and conversations, trad. A. J. Pomerans, Harper \& Row, New York, Evanston and Toronto, 1971.

 \bibitem{Jamblique}  {\sc Jamblique}, 
Vie de Pythagore, Introduction, traduction et notes par L. Brisson \& A. Ph. Segonds, Les Belles Lettres, Paris, 1996. 





\bibitem{Kepler-HM} J. Kepler,  Harmonices Mundi (1619),  traduction fran\c caise de Jean Peyroux,  \emph{L'harmonie du monde}, Librairie A. Blanchard, Paris, 1977.
%
%
 \bibitem{Macrobe} Macrobe, Commentaire au songe de Scipion, Tiré de la République de Cicéron, Paris, Firmin Didot, 1875.

\bibitem{Negrepontis} S. Negrepontis,  Plato on geometry
and the geometers, in S. G. Dani and A. Papadopoulos (eds.), Geometry in history,
Springer Verlag, 2019, p. 1-88.

\bibitem{N2} S. Negrepontis,  The periodic anthyphairetic nature of the One of the Second Hypothesis in  
Plato's Parmenides, Actes du colloque Mathématiques et musique : des Grecs à  Euler, 10-11 September, 2015, Strasbourg, (dir. X. Hascher et A. Papadopoulos), Hermann, Paris, à paraître.

\bibitem{NK} S. Negrepontis and D. Kalisperi, The mystery of Plato's receptacle in the \emph{Timaeus} resolved, in Handbook of the History and Philosophy of Mathematical practice, ed. B. Sriraman, Springer, 2020, p. 1-75.

\bibitem{Nicomaque} Nicomaque de Gérase, Manuel d'harmonique, trad. Ch.-\'E. Ruelle, Annuaire de l'Association pour l'encouragement des \'Etudes grecques en France (1880), Paris, Baur 1881. Trad. anglaise A. Barker, Greek Musical Writings, Cambridge University Press, 1989.

\bibitem{Noll} Th. Noll, Two Notions of well-formedness in the organization of musical pitch, Music\ae \  Scienti\ae , 2010

\bibitem{Pappus} Pappus d'Alexandrie. La collection mathématique ;  traduit du  grec  par Paul Ver Eecke. Paris, Bruges, Desclée de Brouwer, 1933.

\bibitem{Timee} Platon, Timée, In : Platon, \OE uvres complètes, vol. II, trad. L. Robin et M. J. Moreau, Gallimard, Bibliothèque de la Pléiade, Paris, 1950, p. 431-524.

\bibitem{Republique} Platon, La République, In : Platon, \OE uvres complètes, vol. I, trad. L. Robin, et M. J. Moreau, Gallimard, Bibliothèque de la Pléiade, Paris, 1950, p. 858-1241.

\bibitem{Epinomis} Platon, \'Epinomis, In: Platon, \OE uvres complètes, vol. II, trad. L. Robin et M. J. Moreau, Gallimard, Bibliothèque de la Pléiade, Paris, 1950, p. 1133-1163.


\bibitem{Sophiste} Platon, Sophiste, In: Platon, \OE uvres complètes, vol. II, trad. L. Robin et M. J. Moreau, Gallimard, Bibliothèque de la Pléiade, Paris, 1950, p. 257-338.


\bibitem{Lois} Platon, Les Lois, In : Platon, \OE uvres complètes, vol. II, trad. L. Robin et M. J. Moreau, Gallimard, Bibliothèque de la Pléiade, Paris, 1950, p. 634-1131.



 \bibitem{Porphyre} Porphyre, Vie de Pythagore, dir. \'E. Des Places, Paris, Les Belles Lettres, 2003.


 \bibitem{Plutarque} Plutarque, Propos de table (Moralia), éd. et trad. F. Frazier et J. Sirinelli, Paris, Les Belles lettres, 1996.

%
%
%
  \bibitem{Plutarque-Quaestiones} Plutarque, Questions platoniques, dans \cite{Plutarque} livre V.
 

%
%
%
  
%
%
\bibitem{Reinach} Th. Reinach, La musique des sphères, Revue des \'Etudes Grecques, (1900) 13, fascicule 55. p. 432-449.
       
\bibitem{Taylor} A. E. Taylor, A commentary on Platon's Timaeus,  Oxford, Clarendon Press, 1928


 
 \bibitem{Theon}  Th\'eon de Smyrne,   Exposition des connaissances math\'ematiques utiles pour la lecture de Platon,
\'edition  bilingue (grec-fran\c cais) \'etablie par J. Dupuis, Paris 1892, 
reprint par les  \'editions Culture et Civilisation, 115 Av. Gabriel Lebon, Bruxelles
1966.

\bibitem{Vitrac} B. Vitrac, Les mathématiques dans le Timée de Platon : le
point de vue d'un historien des sciences, \'Etudes platoniciennes, 2 (2006) p. 11-78.
%
%
%

\bibitem{Whitehead}  A. N. Whitehead, Process and Reality, corrected edition, ed. D. R. Griffin and D. W. Sherburne, London, New York, 1973.


\end{thebibliography}
\end{document}